\documentclass[12pt]{article}
\usepackage{amssymb,latexsym}
\title{Non-\sa{} perturbations of
\sa{} \op{}s in 2 dimensions II. Vanishing Averages}
\author{Michael Hitrik\\Department of Mathematics \\ University of
California \\Los Angeles \\ CA 90095-1555,
USA\\hitrik@math.ucla.edu \and Johannes Sj\"ostrand\\Centre de
Math\'ematiques\\Laurent Schwartz\\Ecole Polytechnique\\FR--91128
Palaiseau\\France
\\and UMR 7640 CNRS \\johannes@math.polytechnique.fr}
\date{}

\def\wrtext#1{\relax\ifmmode{\leavevmode\hbox{#1}}\else{#1}\fi}
\def\abs#1{\left|#1\right|}
\def\begeq{\begin{equation}}
\def\endeq{\end{equation}}
\def\Remark{\vskip 2mm \noindent {\em Remark}}

\def\ekv#1#2{\begeq\label{#1}#2\endeq}

\def\iint{\int\hskip -2mm\int}

\def\an{analytic}

\def\canform{canonical transformation}

\def\fu{function}

\def\fop{Fourier integral operator}
\def\fourior{Fourier integral operator}
\def\fouriors{Fourier integral operators }
\def\hol{holomorphic}
\def\indep{independent}

\def\mfld{manifold}

\def\neigh{neighborhood}

\def\op{operator}

\def\pseudor{pseudodifferential operator}

\def\sa{selfadjoint}

\def\Re{{\rm Re\,}}
\def\Im{{\rm Im\,}}

\newcommand{\eps}{\epsilon}
\def\part#1{\frac{\partial}{\partial #1}}
\def\half{\frac{1}{2}}

\def\norm#1{||\,#1\,||}

\newcommand{\real}{\mbox{\bf R}}
\newcommand{\comp}{\mbox{\bf C}}
\newcommand{\z}{\mbox{\bf Z}}
\newcommand{\nat}{\mbox{\bf N}}

\renewcommand{\Re}{\mbox{\rm Re\,}}
\renewcommand{\Im}{\mbox{\rm Im\,}}

\renewcommand{\exp}{\mbox{\rm exp\,}}
\newcommand{\supp}{\mbox{\rm supp}}
\newcommand{\T}{\mbox{\bf T}}

\newtheorem{dref}{Definition}[section]

\newtheorem{theo}[dref]{Theorem}
\newtheorem{prop}[dref]{Proposition}

\newenvironment{proof}{\vspace{.3cm}\noindent{{\em Proof:}}}{\hfill$\Box$}

\begin{document}
\maketitle

\begin{abstract} This is the second in a series of works devoted
to small non-selfadjoint perturbations of selfadjoint
semiclassical pseudodifferential operators in dimension 2. As in
our previous work, we consider the case when the classical flow of
the unperturbed part is periodic. Under the assumption that the
flow average of the leading perturbation vanishes identically, we
show how to obtain a complete asymptotic description of the
individual eigenvalues in certain domains in the complex plane,
provided that the strength of the perturbation $\eps$ is $\gg
h^{1/2}$, or sometimes only $\gg h$, and enjoys the upper bound
$\eps={\cal O}(h^{\delta})$, for some $\delta>0$.
\end{abstract}




\section{Introduction}\label{section0}
\setcounter{equation}{0} \setcounter{dref}{0}

This paper is the second one in a series dealing with the
semiclassical spectral asymptotics of small non-selfadjoint
perturbations of selfadjoint operators in dimension 2. In our
previous work~\cite{HiSj} we studied the case when the classical
bicharacteristic flow of the unperturbed operator is periodic in
some energy shell, and in this paper we shall take a next step in
this study, still assuming the periodicity of the underlying
classical flow.

In this entire project we have to a large extent been motivated by
the work of A. Melin and the second author~\cite{MeSj}, where it
has been observed that in dimension 2, there exist wide and stable
classes of non-selfadjoint operators, for which one can obtain a
detailed information about the individual eigenvalues in the
semiclassical limit $h\rightarrow 0$, in some fixed domain in
$\comp$. While the corresponding conclusion is well-known and
classical in the selfadjoint case in dimension 1, important
difficulties appear when passing to higher dimensions, unless one
makes an additional strong assumption of  complete integrability,
also on the quantum level. The latter case has a long tradition in
semiclassical analysis, and we refer to~\cite{CdVVN} for a recent
study of it in dimension 2, in the presence of singularities---see
also the references given in that paper.

The aforementioned difficulties in the selfadjoint case in higher
dimensions are intimately related to the phenomenon of exceptional
sets, corresponding to unstable tori in the KAM theorem, and a
basic discovery of~\cite{MeSj} which allowed for complete spectral
results, was a version of the KAM theorem without small divisors
in the complex domain. Also in~\cite{HiSj}, an important role was
played by certain flow invariant Lagrangian tori in the
complexified phase space, and we have seen how to obtain complete
spectral results in some regions of the complex plane, when the
strength of the non-selfadjoint perturbation $\eps={\cal
O}(h^{\delta})$, $\delta>0$, is bounded from below by an
appropriate power of the semiclassical parameter. Corresponding
lower bounds will play an important role also in this work.

Our general assumptions will be the same as in~\cite{HiSj}, and as
in that paper, we shall let $M$ stand for either $\real^2$ or a
real-analytic compact manifold of dimension 2. We shall also write
$\widetilde{M}$ to denote a complexification of $M$, so that
$\widetilde{M}=\comp^2$ when $M=\real^2$, and in the manifold
case, $\widetilde{M}$ is a Grauert tube of $M$.

When $M=\real^2$, let
\begeq \label{1.1} P_\epsilon
=P(x,hD_x,\epsilon ;h)
\endeq
be the Weyl quantization on ${\bf R}^2$ of a symbol $P(x,\xi
,\epsilon ;h)$ depending smoothly on $\epsilon \in{\rm
neigh\,}(0,{\bf R})$ and taking values in the space of \hol{}
\fu{}s of $(x,\xi )$ in a tubular \neigh{} of ${\bf R}^4$ in ${\bf
C}^4$, with \begeq\label{1.2} \vert P(x,\xi ,\epsilon ;h)\vert \le
Cm(\Re (x,\xi ))
\endeq
there. Here $m$ is assumed to be an order \fu{} on ${\bf R}^4$, in
the sense that $m>0$ and \begeq\label{1.3} m(X)\le C_0\langle
X-Y\rangle ^{N_0}m(Y),\ X,Y\in{\bf R}^4,
\endeq
for some $C_0$, $N_0>0$.  We also assume that \begeq\label{1.4}
m\ge 1.\endeq

Assume furthermore that \begeq\label{1.5} P(x,\xi ,\epsilon
;h)\sim \sum_{j=0}^\infty  p_{j,\epsilon }(x,\xi)h^j,\ h\to 0,
\endeq
in the space of holomorphic \fu{}s. We make the ellipticity
assumption \begeq\label{1.6} \vert p_{0,\epsilon }(x,\xi )\vert
\ge {1\over C}m(\Re (x,\xi )),\ \vert (x,\xi )\vert \ge C,
\endeq
for some $C>0$.

When $M$ is a compact manifold, we let \begeq \label{1.7}
P_{\eps}=\sum_{\abs{\alpha}\leq m}
a_{\alpha,\eps}(x;h)(hD_x)^{\alpha}
\endeq
be a differential operator on $M$, such that for every choice of
local coordinates, centered at some point of $M$,
$a_{\alpha,\eps}(x;h)$ is a smooth function of $\eps$ with values
in the space of bounded holomorphic functions in a complex
\neigh{} of $x=0$. We assume furthermore that
 \begeq
a_{\alpha,\eps}(x;h)\sim \sum_{j=0}^{\infty}
a_{\alpha,\eps,j}(x)h^j,\quad h\rightarrow 0,
\endeq
in the space of such functions. The semiclassical principal symbol
in this case is given by \begeq \label{1.9} p_{0,\eps}(x,\xi)=\sum
a_{\alpha,\eps,0}(x)\xi^{\alpha},
\endeq
and we make an ellipticity assumption
$$
\abs{p_{0,\eps}(x,\xi)}\geq \frac{1}{C}\langle{\xi}\rangle^m,\quad
(x,\xi)\in T^*M,\,\,\abs{\xi}\geq C,
$$
for some large $C>0$. Here we assume that $M$ has been equipped
with some analytic Riemannian metric, so that $\abs{\xi}$ and
$\langle{\xi}\rangle=(1+\abs{\xi}^2)^{1/2}$ are defined.

\par Sometimes, we write $p_\epsilon $ for $p_{0,\epsilon }$ and
simply $p$ for $p_{0,0}$. Assume \begeq\label{1.10} P_{\epsilon
=0} \hbox{ is formally \sa{}.}
\endeq
In the case when $M$ is compact, we let the underlying Hilbert
space be $L^2(M, \mu(dx))$, where $\mu(dx)$ is the Riemannian
volume element.

For $h>0$ small enough and when equipped with the domain $H(m)$,
the natu\-rally defined Sobolev space associated with the weight
$m$ (so that in the compact case, $H(m)$ is the standard Sobolev
space $H^m$ with an $h$ dependent norm), $P_{\eps}$ becomes a
closed densely defined operator on $L^2(M)$. Moreover, the
assumptions above imply that the spectrum of $P_{\eps}$ in a fixed
\neigh{} of $0\in \comp$ is discrete, when $h>0$ and $\eps\geq 0$
are sufficiently small. Clearly, if $z\in {\rm neigh}(0,\comp)$ is
an eigenvalue of $P_{\eps}$, then $\Im z={\cal O}(\eps)$. As
in~\cite{HiSj}, we shall be interested in the asymptotics of
individual eigenvalues of $P_{\eps}$ inside a band $\abs{\Im
z}\leq {\cal O}(\eps)$.

\par Assume for simplicity that (with $p=p_{\epsilon =0}$)
\begeq\label{1.11} p^{-1}(0)\cap T^*M\hbox{ is connected.}\endeq
Let $H_p=p'_\xi \cdot {\partial \over \partial x}-p'_x\cdot
{\partial \over \partial \xi }$ be the Hamilton field of $p$. In
this work we shall always assume that for $E\in{\rm
neigh\,}(0,{\bf R})$:
\begin{eqnarray}\label{1.12}
\hbox{The }H_p\hbox{-flow is periodic on }p^{-1}(E)\cap T^*M\hbox{
with}\\ \hbox{period }T(E)>0 \hbox{ depending \an{}ally on }E.\nonumber
\end{eqnarray}

 When $f$ is an analytic function defined near
$p^{-1}(0)\cap T^*M$, we introduce the trajectory average, \begeq
\label{1.12prime} \langle{f}\rangle=\frac{1}{T(E)}\int_0^{T(E)}
f\circ\exp(tH_p)\,dt\quad \wrtext{on}\quad p^{-1}(E)\cap T^*M,
\endeq
which Poisson commutes with $p$, $H_p\langle{f}\rangle=0$. Let us
write out the first few terms in a Taylor expansion of $p_{\eps}$,
\begeq \label{1.11prime} p_{\eps}=p+i\eps q +\eps^2 r+{\cal
O}(\eps^3m),
\endeq
in the case $M=\real^2$, and $p_{\eps}=p+i\eps q+\eps^2r+{\cal
O}(\eps^3 \langle{\xi}\rangle^m)$ in the compact case.

In~\cite{HiSj}, restricting attention to the range
$$
h\ll \eps={\cal O}(h^{\delta}),
$$
for an arbitrary but fixed $\delta>0$, we have obtained complete
asymptotic expansions for all eigenvalues of $P_{\eps}$ in
rectangles of the form
$$
[-1/{\cal O}(1),1/{\cal O}(1)]+i\eps[F_0-1/{\cal O}(1),F_0+1/{\cal
O}(1)],
$$
under the assumption that $F_0$ is either a regular value of $\Re
\langle{q}\rangle$, viewed as a function on $p^{-1}(0)\cap T^*M$,
or a critical value corresponding to a non-degenerate maximum or
minimum of this function. When the subprincipal symbol of
$P_{\eps=0}$ vani\-shes, we were able to treat even smaller values
of $\eps$, $h^2\ll \eps\leq {\cal O}(h^{\delta})$.

Now there exist many natural situations when the averaged
perturbation $\langle{q}\rangle$ vanishes identically, and in
section 7 of~\cite{HiSj} we have already observed that one such
case occurs when studying barrier top resonances for the
semiclassical Schr\"odinger operator, in the case of a $1:1$
resonance for the fundamental frequencies at the top of the
barrier. The purpose of the present work is to carry out a
spectral analysis of $P_{\eps}$ under the assumption that
$\langle{q}\rangle\equiv 0$. This case has a long tradition in the
study of selfadjoint operators with a periodic classical flow
(\cite{Gu1},~\cite{Gu2},~\cite{Iv},~\cite{Uribe}), and at this
point we would also like to mention an interesting paper by L.
Friedlander~\cite{Fr}, which studies the spectral localization for
the Laplace operator on the $n$-sphere, perturbed by a
complex-valued bounded odd potential---see~\cite{Gu1} for a
corresponding study in the real case. A common feature of all the
works mentioned above (although not of~\cite{Fr}), which was also
exploited in~\cite{HiSj}, is an averaging procedure, which allows
one to reduce the dimension by one unit. (For that it is important
that we have $\eps={\cal O}(h^{\delta})$, $\delta>0$.) Here we
also follow this idea, although, as we shall see, some
modifications in the averaging procedure, as compared
with~\cite{HiSj}, will be required. We would also like to
emphasize that the implementation of the averaging method in the
present non-selfadjoint setting involves conjugating our operator
by means of Fourier integ\-ral operators with complex phase,
suitably realized on the FBI transform side, and as
in~\cite{HiSj}, we are naturally led to work in modified
exponentially weighted spaces, in the spirit
of~\cite{Sj6},~\cite{HeSj},~\cite{MeSj}. Additional modifications
of the spaces will be required, and for the final spectral
results, we shall sometimes also have to introduce a supplementary
global hypothesis of dynamical nature, which will ultimately allow
the spectral problem to become microlocalized to a small \neigh{}
of a suitable Lagrangian torus. Let us also remark that we expect
that global hypotheses of a similar kind will appear naturally in
further works in our series, in more general situations, when the
energy surface of the unperturbed symbol admits certain invariant
Lagrangian tori with good arithmetic (diophantine) properties. At
many essential points however, the proofs will be the same as
in~\cite{HiSj}, and for that reason, the presentation of those
parts of the proofs will be somewhat less detailed.

The plan of the paper is as follows.

\noindent In section 2, we reexamine the averaging procedure
of~\cite{HiSj}, adapting it to the present case, and describe a
localization of the spectrum of $P_{\eps}$.

\noindent In section 3 we construct a quantum Birkhoff normal form
near a suitable torus, and, provided that a certain averaged
correction has a nondegenerate imaginary part, obtain the first
spectral result. This is done completely along the lines of the
analysis of~\cite{HiSj}.

\noindent In section 4, we study the complementary case when the
relevant averaged correction is real on the real domain, and we
then show how to obtain complete spectral results, by introducing
an additional global dynamical assumption.

\noindent In section 5 we apply the results of section 3 to
barrier top resonances, thereby complementing the discussion of
section 7 of~\cite{HiSj}.

\noindent In section 6, we give a first application of the results
of section 4 to the spectrum of the damped wave equation on the
2-sphere (~\cite{AsLe},~\cite{Le},~\cite{Sj3},~\cite{Hi}).

\vskip 4mm \noindent {\bf Acknowledgment}. The first author
gratefully acknowledges the support of the MSRI postdoctoral
fellowship as well as a partial support of the National Science
Foundation under the grant DMS--0304970.

\section{Averaging reduction and spectral
localization}\label{section1} \setcounter{equation}{0} In this
section we shall work near $p^{-1}(0)\cap T^*M$ and the following
arguments will be valid in any dimension $n\geq 2$.

Let us recall that in a \neigh{} of $p^{-1}(0)\cap T^*M$, the
principal symbol of $P_{\eps}$ has the form
$$
p_{\eps}=p+i\eps q+\eps^2 r+{\cal O}(\eps^3).
$$
Assuming that $\langle{q}\rangle\equiv 0$, we shall first
reexamine the reduction by averaging, described in section 3
of~\cite{HiSj}, and adapt it to the present situation. Let $G_0$
be an analytic function defined near $p^{-1}(0)\cap T^*M$ such
that \begeq \label{2.1} H_p G_0=q.
\endeq

As remarked in~\cite{HiSj}, we may take \begeq \label{2.2}
G_0=\frac{1}{T(E)}\int_0^{T(E)} t
\left(q\circ\exp(tH_p)\right)\,dt \quad \wrtext{on}\quad
p^{-1}(E)\cap T^*M.
\endeq

Replace now $T^*M$ by the new IR-\mfld{} $\Lambda_{\eps
  G_0}=\exp(i\eps H_{G_0})(T^*M)$, which is
defined in a complex \neigh{} of $p^{-1}(0)\cap T^*M$. Writing
$$
(x,\xi)=\exp(i\eps H_{G_0})(y,\eta),
$$
and using
$\rho=(y,\eta)$ as real symplectic coordinates on
$\Lambda_{\eps G_0}$, we get
\begin{eqnarray}
\label{2.3} \left(p_{\eps}\right)|_{\Lambda_{\eps G_0}}
& = &
\left(p+i\eps q+\eps^2 r+{\cal
O}(\eps^3)\right)\left(\exp(i\eps H_{G_0})(\rho)\right) \\
\nonumber & = & \sum_{k=0}^{\infty} \frac{(i\eps
H_{G_0})^k}{k!}(p+i\eps q+\eps^2 r+{\cal O}(\eps^3)) \\ \nonumber
& = & p+\eps^2\left(r-H_{G_0} q-\half H_{G_0}^2 p\right)+{\cal
O}(\eps^3).
\end{eqnarray}
We shall assume in what follows that $\eps ={\cal O}(h^{\delta})$ for some
fixed $\delta>0$, and iterating the procedure above, we shall modify
the choice of the weight function $G_0$ and of the corresponding
IR-deformation, in order to improve (\ref{2.3}). Set
\begeq
\label{2.4}
G=G_0+i\eps G_1+{\cal O}(\eps^2),
\endeq
where $G_1$ is to be chosen, and compute the restriction of $p_{\eps}$ to the corresponding IR-\mfld{}
$\Lambda_{\eps G}$. We get
\begin{eqnarray*}
p_{\eps}(\exp(i\eps H_G))  & = & p+i\eps q+\eps^2 r+i\eps
H_G(p+i\eps q)-\frac{\eps^2}{2} H_G^2 p +
{\cal O}(\eps^3) \\
& = & p+i\eps q+\eps^2 r+i\eps\left(H_{G_0}+i\eps
H_{G_1}\right)(p+i\eps q)-\frac{\eps^2}{2}H_{G_0}^2 p +{\cal
O}(\eps^3) \\
& = & p+i\eps q +i\eps\left(H_{G_0}p+i\eps H_{G_0}q+i\eps H_{G_1}p\right)+\eps^2 r-\frac{\eps^2}{2} H^2_{G_0}p+{\cal O}(\eps^3) \\
& = & p + \eps^2\left(r+H_p G_1-H_{G_0} q - \half H^2_{G_0}p\right)+{\cal O}(\eps^3).
\end{eqnarray*}
Using that $H_p G_0=q$, we see that the expression in front of $\eps^2$ takes the form
$$
s:=r+H_p G_1-\half H_{G_0}q.
$$
We choose $G_1$ as an analytic solution in a \neigh{} of
$p^{-1}(0)\cap T^*M$ of the equation
$$
H_p G_1=\half\left(H_{G_0}q-\langle{H_{G_0}q\rangle}\right)-(r-\langle{r}\rangle),
$$
and then get a reduction of the original symbol to
\begeq
\label{2.5}
p+\eps^2\langle{s}\rangle +{\cal O}(\eps^3),
\endeq
where
$$
\langle{s}\rangle=\langle{r}\rangle-\half\langle{H_{G_0}
  q\rangle}=\langle{r}\rangle-\frac{1}{2T(E)}\int_0^{T(E)}
  \{G_0,q\}\circ \exp(tH_p)\,dt,
$$
so that $p$ and $\langle{s}\rangle$ are in involution. Recalling
the expression for $G_0$ from (\ref{2.2}) and writing $T(p)$
rather than $T(E)$ for the period on the energy surface
$p^{-1}(E)$, we get
\begin{eqnarray*}
& & \{G_0,q\}=\frac{1}{T(p)}\int_0^{T(p)} s\{q\circ
\exp(sH_p),q\}\,ds\\
& + & \{T(p),q\} \left(q-\frac{1}{T^2(p)}\int_0^{T(p)} s q\circ
\exp(sH_p)\,ds\right).
\end{eqnarray*}
It follows that
\begin{eqnarray}\label{2.6}
& & \langle{s}\rangle=\langle{r}\rangle
\\ \nonumber
& - & \frac{1}{2T(p)^2}\int\!\!\!\int
u\{q\circ\exp((u+v)H_p),q\circ\exp(vH_p)\}\,du\,dv \\ \nonumber &-
& \frac{1}{2T(p)}\int_0^{T(p)} q\circ \exp(v H_p)\{T(p),q\circ
\exp(vH_p)\}\,dv \\ \nonumber & + & \frac{1}{2T^3(p)}
\int\!\!\!\int u\{T(p),q\circ \exp(vH_p)\} q\circ
\exp((u+v)H_p)\,du\,dv,
\end{eqnarray}
where the integration in the double integrals in the right hand
side of (\ref{2.6}) is performed over the rectangle $[0,T(p)]^2$.
Notice also that any other choice of $G_0$ satisfying (\ref{2.1})
gives the same expression for $\langle{s}\rangle$. Indeed, if $h$
is a function invariant under the $H_p$--flow, then
\begin{eqnarray*}
\langle{\{h,q\}}\rangle & = & \frac{1}{T(p)}\int_0^{T(p)}
\{h,q\}\circ \exp(tH_p)\,dt \\
& = & \frac{1}{T(p)}\int_0^{T(p)}\{h,q\circ
\exp(tH_p)\}\,dt=\{h,\langle{q}\rangle\}=0.
\end{eqnarray*}

For future reference, we remark that since $\langle{q}\rangle=0$,
it is true that
$\langle{G_0}\rangle=0$. Indeed,
$$
\langle{G_0}\rangle=\frac{1}{T(E)^2} \int_0^{T(E)} u\left(\int_0^{T(E)} q\circ\exp((u+v)H_p)\,dv\right)\,du=0.
$$

After replacing $p_{\eps}$ by $p_{\eps}\circ \exp(i\eps H_G)$, and
correspondingly $P_{\eps}$ by $U_{\eps}^{-1}P_{\eps}U_{\eps}$,
where $U_{\eps}=e^{\eps G(x,hD_x)/h}$ is a \fourior{} quantizing
$\exp(i\eps H_{G})$, which is defined microlocally near
$p^{-1}(0)\cap T^*M$, we obtain a reduction of our operator
$P_{\eps}$ to an operator with the leading symbol
$$
p_{\eps}=p+\eps^2 \langle{s}\rangle+{\cal O}(\eps^3).
$$
Moreover, if the subprincipal symbol of $P_{\eps}$ is ${\cal
O}(\eps)$, then this is also true after the conjugation, in view of
the improved Egorov property of $U_{\eps}$---see section 2
of~\cite{HiSj}.

As in section 3 of~\cite{HiSj}, let $g:{\rm
neigh}(0,\real)\rightarrow \real$ be the analytic function defined
by \begeq
\label{2.6.5} g'(E)=\frac{T(E)}{2\pi},\quad g(0)=0,
\endeq
so that $g\circ p$ has a $2\pi$-periodic Hamilton flow. Set
$f=g^{-1}$.
\begin{prop}
Assume that the subprincipal symbol of $P_{\eps=0}$ vanishes and that
$\eps\ll h^{1/2}$. Then
the spectrum of $P_{\eps}$ near $0$ is contained in the union of the rectangles of the form
\begin{eqnarray}
\label{2.7}
& &
I_k(\eps)=f\left(h\left(k-\frac{\alpha}{4}\right)-\frac{S}{2\pi}\right)+\left[-{\cal
O}(\eps^2+h^2),{\cal O}(\eps^2+h^2)\right]\\ \nonumber
& + & i[-{\cal O}(\eps^2+\eps h), {\cal O}(\eps^2+\eps h)],
\end{eqnarray}
for $k\in \z$. Here $\alpha\in \z$ and $S\in \real$ are the Maslov
index and classical action, respectively, computed along a closed
$H_p$-trajectory $\subset p^{-1}(0)\cap T^*M$.
\end{prop}
\Remark. The vanishing of the subprincipal symbol of $P_{\eps=0}$
in this result could be weakened to assuming that the trajectory
average of it is constant on each energy surface $p^{-1}(E)\cap
T^*M$, $E\in {\rm neigh}(0,\real)$. Notice also that since $\eps
\ll h^{1/2}$ and $f'(0)>0$, it is true that $I_k(\eps)\cap
I_{k'}(\eps)=\emptyset$, when $k\neq k'$ and $h$ is small enough.

\begin{proof}
The averaging reduction described above together with the spectral
theorem shows that it suffices to prove the proposition when
$\eps=0$, in which case it is well-known and was established in
the semiclassical case in~\cite{HeRo}, following the earlier
works~\cite{We} and~\cite{Co}. Here we shall take the opportunity
to sketch a proof of this result, which does not rely upon the
averaging procedure on the operator level as in~\cite{HeRo}, but
rather works directly with microlocal representatives near the
closed orbits. Similar arguments will be used in what follows. In
doing so, we shall have to assume that $T(0)$ is the minimal
period of every closed $H_p$--trajectory in $p^{-1}(0)\cap T^*M$.
It is quite likely however, that the following argument can be
modified to cover the general case when there exist subperiodic
orbits, and we plan to return to that in the future.

Let $\gamma\subset p^{-1}(0)\cap T^*M$ be a closed
$H_p$-trajectory of minimal period $T(0)$. From section 3
of~\cite{HiSj} we recall that there exists a unitary \fourior{},
microlocally defined from a \neigh{} of $\gamma$ in $T^*M$ to a
\neigh{} of $\tau=x=\xi=0$ in $T^*(S^1_t\times \real_x^{n-1})$,
such that after the conjugation by this operator, the operator
$P=P_{\eps=0}$ takes the form $\widetilde{P}=f(hD_t)+{\cal
O}(h^2)$, where the remainder is an $h$-\pseudor{} whose symbol is
${\cal O}(h^2)$. Here we also use that the Fourier integral
operator has the improved Egorov property. The operator $\widetilde{P}$ acts on the space
$L^2_S$ of functions defined microlocally near $\tau=x=\xi=0$ in
$T^*(S^1_t\times \real_x^{n-1})$, and satisfying the Floquet-Bloch
condition
$$
u(t-2\pi,x)=e^{i(\frac{S}{h}+\frac{2\pi \alpha}{4})}u.
$$
If now $z\in {\rm neigh}(0,\real)$ avoids the
union of the intervals $I_k(0)$, then we see that the operator
\begeq
\label{0.30}
\widetilde{P}-z=f(hD_t)+{\cal O}(h^2)-z
\endeq
is invertible, microlocally near $\tau=x=\xi=0$, with the norm of the
inverse being ${\cal O}(h^{-2})$.

Take now finitely many closed trajectories $\gamma_1,\ldots
\gamma_N \subset p^{-1}(0)\cap T^*M$ and small open
$H_p$--invariant neighborhoods $\Omega_j$ of $\gamma_j$, $1\leq
j\leq N$, such that the $\Omega_j$'s form a cover of
$p^{-1}(0)\cap T^*M$. Take also functions $\chi_j\in
C^{\infty}_0(\Omega_j)$, $1\leq j\leq N$, such that $\sum_{j=1}^N
\chi_j=1$ near $p^{-1}(0)\cap T^*M$, and the $\chi_j$ are in
involution with $p$. When $z\in {\rm neigh}(0,\real)$ avoids the
intervals $I_k(0)$, we consider the equation
$$
(P-z)u=v,\quad u\in H(m).
$$
In what follows we write $U_j$ to denote the \fourior{} intertwining
$P$ and the operator $\widetilde{P}$, microlocally in
$\Omega_j$. Then
$$
(\widetilde{P}-z)U_j\chi_j u=U_j\left(\chi_j v+[P,\chi_j]u\right),
$$
modulo an ${\cal O}(h^{\infty})$-error. Here $[P,\chi_j]={\cal
O}(h^3)$ in the operator sense, since the subprincipal symbols of $P$
and $\chi_j$ both vanish. It follows that, with $\norm{\cdot}$
denoting the $L^2$--norm,
$$
\norm{\chi_j u}\leq {\cal O}\left(\frac{1}{h^2}\right)\norm{v}+{\cal O}(h)\norm{u}.
$$
Summing these estimates over all $j$ and combining them with the
usual elliptic bound away from $p^{-1}(0)\cap T^*M$, we conclude
that
$$
\norm{u}\leq {\cal O}\left(\frac{1}{h^2}\right)\norm{v},
$$
when $h$ is small enough. The injectivity and hence the
invertibility of $P-z$, for $z\in {\rm neigh}(0,\real)$ avoiding
the rectangles $I_k(0)$, $k\in \z$, follows.
\end{proof}

\Remark. If $\langle{s}\rangle$ is real on the real domain, then we
obtain an improvement on the vertical size of the rectangles
$I_k(\eps)$ in (\ref{2.7}), and we see that if $z\in {\rm Spec}(P_{\eps})\cap
{\rm neigh}(0,\comp)$, then
$$
\Im z={\cal O}(\eps^3+\eps h),
$$
for $\eps={\cal O}(h^{\delta})$, $\delta>0$.

\section{Normal forms and eigenvalues in the torus case I}\label{section3}
\setcounter{equation}{0} The arguments in this section will to a
large extent be parallel to the corresponding discussion in
sections 4 and 6 of~\cite{HiSj}, and therefore the following
presentation will be less detailed.

We have seen in section 2 that we may reduce ourselves to the case
of an operator with the leading symbol
$$
p_{\eps}=p+\eps^2\langle{s}\rangle+{\cal O}(\eps^3).
$$
Throughout this section, we shall assume that
\begin{eqnarray}
\label{3.0.5} & & \Re\langle{s}\rangle \, \wrtext{is an analytic
function of}\, \Im \langle{s}\rangle \, \wrtext{and}\, p, \\
\nonumber
 & & \wrtext{in the region of $T^*M$ where}\, \abs{p}\leq
1/{\cal O}(1).
\end{eqnarray}

Introduce the $H_p$-invariant set
$$
\Lambda_{0,0}: p=0,\, \Im \langle{s}\rangle=0,
$$
and assume that the minimal period for all the closed
$H_p$--trajectories in $\Lambda_{0,0}$ is $T(0)$, and that $dp$
and $d \langle{\Im s}\rangle$ are linearly independent at each
point of $\Lambda_{0,0}$. Then $\Lambda_{0,0}$ is a Lagrangian
manifold which is also a union of tori. Assume for simplicity that
$\Lambda_{0,0}$ is connected, so that it is equal to a single
Lagrangian torus. In a \neigh{} of $\Lambda_{0,0}$, $p$ and
$\Im\langle{s}\rangle$ form a completely integrable system, and
precisely as in section 4 of~\cite{HiSj}, we perform a real
analytic canonical transformation \begeq \label{3.1} \kappa: {\rm
neigh}\left(\xi=0, T^*{\bf T}^2\right)\rightarrow {\rm
neigh}\left(\Lambda_{0,0},T^*M\right),
\endeq
such that $p\circ\kappa=p(\xi_1)$, $\Im
\langle{s}\rangle\circ\kappa=\Im \langle{s}\rangle(\xi)$. (In
section 4 of~\cite{HiSj} we recalled the construction of
$\kappa$---see also~\cite{Du}.) It is therefore clear that when
working microlocally near $\Lambda_{0,0}$ and implementing
$\kappa$ by means of a microlocally unitary \fourior{}, which also
has the improved Egorov property, we reduce ourselves to a new
operator, still denoted by $P_{\eps}$, which is microlocally
defined near the zero section in $T^*{\bf T}^2$, and which has the
principal symbol
$$
p(\xi_1)+\eps^2 \langle{s}\rangle(\xi)+{\cal O}(\eps^3).
$$
Also, the transformed operator will act on the space
$L^2_{\theta}(\T^2)$ of microlocally defined Floquet periodic
functions on ${\bf T}^2$, satisfying
$$
u(x-\nu)=e^{i\theta\cdot \nu} u(x),\quad \nu \in (2\pi \z)^2,\quad
\theta=\frac{S}{2\pi h}+\frac{\alpha}{4}.
$$
Here $S=(S_1,S_2)$ with $S_j$ being the action associated to the
fundamental cycle $\gamma_j$ in $\Lambda_{0,0}$, $j=1,2$, with $\gamma_1$
being given by a closed $H_p$--trajectory, and
$\alpha=(\alpha_1,\alpha_2)$ is the corresponding Maslov
index. The assumption about the linear independence of the differentials implies
$$
\partial_{\xi_1}p(0)\neq 0, \quad \partial_{\xi_2}\Im \langle{s}\rangle(0)\neq 0.
$$

Using the first of these assumptions and repeating the arguments
of section 4 of~\cite{HiSj}, we find an $h$-\pseudor{}
$A=\sum_{\nu=0}^{\infty} h^{\nu} a_{\nu}(x,\xi,\eps)$, with
$a_0={\cal O}(\eps^3)$, such that the full symbol of the operator
$$
\widetilde{P}=e^{\frac{i}{h} A}P_{\eps} e^{-\frac{i}{h}
  A}=e^{\frac{i}{h}{\rm ad} A}P_{\eps}
=\sum_{k=0}^{\infty} \frac{1}{k!}
\left(\frac{i}{h}{\rm ad}\,A\right)^k P_{\eps}
$$
is independent of $x_1$. We remark that the functions
$a_{\nu}(x,\xi,\eps)$ are constructed as formal power series in
$\eps$, with coefficients holomorphic in a fixed complex \neigh{}
of $\xi=0$. These formal power series are then realized as
$C^{\infty}$-symbols since we work under the assumption that
$\eps={\cal O}(h^{\delta})$ for some fixed $\delta>0$. We are thus
reduced to the operator
\begeq \label{3.2}
\widetilde{P}_{\eps}=\sum_{j=0}^{\infty} h^j
\widetilde{p}_j(x_2,\xi,\eps),
\endeq
with
$$
\widetilde{p}_0=p(\xi_1)+\eps^2 \langle{s}\rangle(\xi)+{\cal
O}(\eps^3),
$$
also independent of $x_1$. We would like to perform further
conjugations of $\widetilde{P}_{\eps}$ by means of Fourier
integral operators, in order to make its leading symbol
independent of $x_2$ as well. We shall do so first in the case
when the subprincipal symbol of $P_{\eps=0}$ is not necessarily
zero. On the symbol level, we write
\begin{eqnarray}
\label{3.3} \widetilde{P}_\epsilon  & = & p(\xi _1)+ \epsilon^2
(\langle s\rangle (\xi )+{\cal O}(\epsilon )+{h\over
\epsilon^2}\widetilde{p}_1(x_2,\xi ,\epsilon )+h{h\over
\epsilon^2}\widetilde{p}_2(x_2,\xi ,\epsilon
)+\ldots )\\
& = & p(\xi _1)+\epsilon^2 (r_0(x_2,\xi ,\epsilon ,{h\over
\epsilon^2})+hr_1(x_2,\xi ,\epsilon ,{h\over \epsilon^2
})+\ldots),\nonumber\end{eqnarray}
with
$$r_0(x_2,\xi ,\epsilon ,{h\over \epsilon^2})=\langle s\rangle
(\xi )+{\cal O}(\epsilon )+{h\over \epsilon^2}\widetilde{p}_1=\langle
s\rangle (\xi )+{\cal O}(\epsilon )+{\cal O}({h\over \epsilon^2
}),$$
$$r_1={h\over \epsilon^2}\widetilde{p}_2(x_2,\xi ,\epsilon ),\, \ldots $$
Notice that $r_j={\cal O}(h/\epsilon^2)$ for $j\ge 1$.

Following~\cite{HiSj}, we shall
treat $h/\epsilon^2$ as an independent small parameter. In the case when
$h/\eps^2 \leq h^{\delta_1}$, $\delta_1>0$, an inspection of the arguments of
section 4 of~\cite{HiSj} shows that there exist operators
$$
B_0=b_0(x_2,hD_x,\epsilon ,h/\epsilon^2),\ b_0={\cal O}(\epsilon
+h/\epsilon^2),
$$ and
$$
B_1=\sum_{\nu =1}^\infty  b_{\nu}(x_2,hD_x,\epsilon ,h/\epsilon^2
)h^{\nu},\ b_{\nu}={\cal O}(\epsilon +h/\epsilon^2),
$$
such that
\ekv{3.4}
{
\widehat{P}_\epsilon :=e^{{i\over h}{\rm ad}_{B_1}}
e^{{i\over h}{\rm ad}_{B_0}}\widetilde{P}_\epsilon }
has a symbol \indep{} of $x$:
\ekv{3.5}{\widehat{P}_\epsilon =p(\xi _1)+\epsilon^2 (r_0(\xi
,\epsilon ,{h\over \epsilon^2})+hr_1(\xi ,\epsilon ,{h\over
\epsilon^2})+\ldots),} with
$r_0=\langle s\rangle (\xi )+{\cal O}(\epsilon +h/\epsilon^2 )$,
and $r_\nu ={\cal O}(\epsilon +h/\epsilon^2)$ for
$\nu \ge 1$.

Remaining in the general case, we shall now discuss  the
construction of the conjugating operators assuming merely that
$h/\eps^2\leq \delta_0$, for some $\delta_0>0$ small enough but
fixed. Again, this is done completely along the lines of the
analysis of~\cite{HiSj}. We therefore find that there exists a
\hol{} \canform{} $\widetilde{\kappa}$ from a complex \neigh{} of
the zero section in $T^*\T^2$ to a similar set, with a generating
function of the form
\begeq \label{3.6} \psi \left(x,\eta
,\epsilon ,{h\over \epsilon^2}\right)=x\cdot \eta +\psi _{\rm
per}\left(x_2,\eta,\eps,\frac{h}{\eps^2}\right) , \quad \psi_{\rm
per}={\cal O}\left(\eps+h/\eps^2\right),
\endeq
such that
\begeq
\label{3.7}
(r_0\circ \kappa )(y,\eta
,\epsilon ,{h\over \epsilon^2})= \langle r_0(\cdot ,\xi ,\epsilon
,{h\over \epsilon^2})\rangle =\langle r_0(\cdot ,\eta ,\epsilon
,{h\over \epsilon^2})\rangle +{\cal O}(\epsilon ^2+
({h\over \epsilon^2
})^2)
\endeq
is independent of $y$. Here the average of $r_0$ is with respect to
the $x_2$ variable.

\par We can quantize $\kappa $ as a \fop{} $U$ and after
conjugation by this \op{}, we may assume that we have a new \op{}
$\widetilde{P}_\epsilon $ as in (\ref{3.3}), with
$r_0=\langle s\rangle (\xi )+{\cal O}(\epsilon +h/\epsilon^2)$
\indep{} of $x$ and with $r_j={\cal O}(\epsilon +h/\epsilon^2)$,
$j\geq 1$.

\par As before, we can then make a further conjugation $e^{{i\over
h}{\rm ad}_{B_1}}$ in order to remove the $x$-dependence
completely, and the conclusion is that if we make no assumption on
the subprincipal symbol and restrict the attention to
$h/\epsilon^2\le \delta _0$, for $\delta _0>0$ small enough, then
we can find a  \fop{}, \ekv{3.8} {U^{-1}u(x;h)={1\over (2\pi
h)^2}\iint e^{{i\over h}(\psi (x,\eta )-y\cdot \eta )} a(x,\eta;h)
u(y)dyd\eta ,} with $\psi (x,\eta )=x\cdot \eta +\psi _{\rm
per}(x_2,\eta ,\epsilon ,h/\epsilon^2)$, $\psi _{\rm per}={\cal
O}(\epsilon +h/\epsilon^2)$, and
$$B_1=\sum_{\nu =1}^\infty  b_{\nu}(x_2,hD_x,\epsilon ,{h\over
\epsilon^2})h^{\nu},\ b_{\nu}={\cal O}(\epsilon +{h\over \epsilon^2}),$$
such that
$$\widehat{P}_\epsilon :=e^{{i\over h}{\rm ad}_{B_1}}{\rm Ad}_U
\widetilde{P}_\epsilon $$
has a symbol \indep{} of $x$ as in (\ref{3.5}), with the same estimates
as there.

\vskip 2mm
\noindent
We shall now turn the attention to the case when the subprincipal
symbol of $P_{\eps=0}$ vanishes. The arguments of~\cite{HiSj} together with
the preceding discussion make it clear that in this case we may reach the range
$h/\eps\leq \delta_0$, for $\delta_0>0$ small enough. Indeed, we write
as in (\ref{3.3}), using that $\widetilde{p}_1=\eps q_1$,
\begin{eqnarray*}
\widetilde{P}_\epsilon & = & p(\xi _1)+\epsilon^2 (\langle s\rangle
(\xi )+{\cal O}(\epsilon )+\frac{h}{\eps} q_1(x_2,\xi ,\epsilon )+{h^2\over
\epsilon^2 }\widetilde{p}_2+h{h^2\over \epsilon^2 }\widetilde{p}_3+\ldots)\\
& = & p(\xi _1)+\epsilon^2 (r_0(x_2,\xi ,\epsilon ,{h\over \epsilon
})+hr_1(x_2,\xi ,\epsilon ,{h\over \epsilon })+h^2r_2+\ldots),
\nonumber\end{eqnarray*}
with
\begin{eqnarray*}
r_0(x_2,\xi ,\epsilon ,{h\over \epsilon }) & = & \langle s\rangle
(\xi )+{\cal O}(\epsilon )+\frac{h}{\eps} q_1+\frac{h^2}{\epsilon^2}\widetilde{p}_2, \\
r_1(x_2,\xi ,\epsilon ,{h\over \epsilon })& = & {h^2\over \epsilon^2}\widetilde{p}_3, \\
 r_2(x_2,\xi ,\epsilon ,{h\over\epsilon })& = & {h^2\over \epsilon^2
}\widetilde{p}_4, \ldots
\end{eqnarray*}
It follows that under the stated smallness assumption on $\eps$,
we can eliminate the $x_2$--dependence  by means of conjugations
by \fouriors{}. To be precise, as in~\cite{HiSj}, we find that
there exist an elliptic \fop{} $U^{-1}$ and $B_1(x_2,hD_x,\epsilon
,h/\epsilon ;h)$ such that
$$
e^{{i\over h}{\rm ad}_{B_1}}{\rm }Ad_{U} \widetilde{P}_\epsilon
=\widehat{P}(hD_x,\epsilon ,h/\epsilon ;h)
$$
has a symbol $\widehat{P}(\xi ,\epsilon ,h/\epsilon
;h)$ of the form
$$
p(\xi_1)+\eps^2 \left(r_0+hr_1+\ldots\right),
$$
where $r_0=\langle{s}\rangle+{\cal O}(\eps)+{\cal O}(h/\eps)$, and
$r_j={\cal O}(\eps+h/\eps)$, $j\geq 1$. This quantum microlocal Birkhoff normal
form near the torus $\Lambda_{0,0}$ leads to the formal
quasi-eigenvalues associated with $\Lambda_{0,0}$ and given by
$\widehat{P}(h(k-\alpha/4)-S/2\pi,\eps,h/\eps~;h)$, $k\in \z^2$.

\vskip 2mm \noindent So far our discussion has been completely
analogous to the corresponding one in section 4 of~\cite{HiSj},
and it is also clear that the analysis of the globally well-posed
Grushin problem from section 6 of that paper can be applied
without any change. To fix the ideas, in what follows, let us
discuss the spectral results assuming that the subprincipal symbol
of $P_{\eps=0}$ vanishes and $h/\eps\leq \delta_0 \ll 1$,
$\eps={\cal O}(h^{\delta})$, $\delta>0$. We then inspect all the
steps of the argument of section 6 of~\cite{HiSj} to find a new
($h$-dependent) globally defined Hilbert space
$H(\widehat{\Lambda}_{\eps})$, associated with a suitable
IR-manifold $\widehat{\Lambda}_{\eps}$, $(\eps+h/\eps)$--close to
$T^*M$, such that the action of $P_{\eps}$ on
$H(\widehat{\Lambda}_{\eps})$ can be, microlocally near
$\Lambda_{0,0}$, identified with the action of
$\widehat{P}(hD_x,\eps,h/\eps;h)$ on $L^2_{\theta}({\bf T}^2)$, by
means of an elliptic semiclassical \fourior{}
$$
U={\cal O}(1): H(\widehat{\Lambda}_{\eps})\rightarrow
L^2_{\theta}(\T^2).
$$
We refer to Proposition 6.1 of~\cite{HiSj} for the precise
description of the basic properties of $U$. As in~\cite{HiSj}, the
construction of the deformation $\widehat{\Lambda}_{\eps}$ is done
in such a way that for $\rho\in \widehat{\Lambda}_{\eps}$ away
from a small \neigh{} of $\Lambda_{0,0}$,
\begeq \label{3.9}
\abs{\Re P_{\eps}(\rho,h)}\geq \frac{1}{{\cal O}(1)}\quad
\wrtext{or}\quad \abs{\Im P_{\eps}(\rho,h)}\geq
\frac{\eps^2}{{\cal O}(1)}.
\endeq

Let now $z$ vary in the rectangle \begeq \label{3.10} \abs{\Re
z}<\frac{1}{C},\quad \abs{\Im z}< \frac{\eps^2}{C},
\endeq
for a sufficiently large constant $C>0$. If $z$ in (\ref{3.10}) is
away from any $\eps^2 h/\widetilde{C}$-neighborhood of the values
$\widehat{P}(h(k-\alpha/4)-S/2\pi, \eps,h/\eps;h)$, $k\in \z^2$,
it follows from the argument explained in section 6 of~\cite{HiSj}
that the operator
$$
P_{\eps}-z: H(\widehat{\Lambda}_{\eps})\rightarrow H(\widehat{\Lambda}_{\eps})
$$
is invertible, with the norm of the inverse ${\cal O}((\eps^2
h)^{-1})$. The setup of the Grushin problem is also identical to
the one in~\cite{HiSj}, and repeating the arguments of that paper
we find that up to ${\cal O}(h^{\infty})$, the eigenvalues of
$P_{\eps}$ in (\ref{3.10}) are given by the quasi-eigenvalues
associated to $\Lambda_{0,0}$,
$\widehat{P}(h(k-\frac{\alpha}{4})-\frac{S}{2\pi},\eps,\frac{h}{\eps};h)$,
$k\in \z^2$. The corresponding result is also true in the range
$h^{1/2}\ll \eps ={\cal O}(h^{\delta})$, $\delta>0$, if the
subprincipal symbol of $P_{\eps=0}$ is not necessarily zero.

We summarize the results of this section in the following theorem.

\begin{theo} Let $P_{\eps}$ be an operator as in the introduction, so that {\rm (\ref{1.12})}, {\rm (\ref{1.11prime})}
hold. Assume that the trajectory average $\langle{q}\rangle$
vanishes. Define $s$ as in section {\rm 2}, so that
$\langle{s}\rangle$ is given by {\rm (\ref{2.6})}, and assume {\rm
(\ref{3.0.5})}. Let $F_0$ be a regular value of $\Im
\langle{s}\rangle$, considered as a function on $p^{-1}(0)\cap
T^*M$. Assume that the Lagrangian manifold
$$
\Lambda_{0,F_0}: p=0,\, \Im \langle{s}\rangle=F_0
$$
is connected, so that it is diffeomorphic to a single torus, and
that $T(0)$ is the minimal period for the $H_p$--flow in
$\Lambda_{0,F_0}$. We fix a basis for the first homology group of
$\Lambda_{0,F_0}$ given by the cycles $\gamma_j$, $j=1,2$, with
$\gamma_1$ corresponding to a closed $H_p$--orbit. Let us write
$S=(S_1,S_2)$ and $\alpha=(\alpha_1,\alpha_2)$ for the actions and
Maslov indices of the cycles, respectively. Assume that
$h^{1/2}\ll \eps={\cal O}(h^{\delta})$ for some $\delta>0$, and
let $C>0$ be large enough. A complex number $z$ satisfying

$$
\abs{\Re z}<\frac{1}{C},\quad \abs{\Im z}< \frac{\eps^2}{C},
$$
is in the spectrum of $P_{\eps}$ precisely when
$$
z=z(k,\eps;h)=\widehat{P}\left(h\left(k-\frac{\alpha}{4}\right)-\frac{S}{2\pi},\eps,\frac{h}{\eps^2};h\right)+{\cal
O}(h^{\infty}), \quad k\in \z^2.
$$
The function $\widehat{P}(\xi,\eps,\frac{h}{\eps^2};h)$ is
holomorphic in $\xi\in {\rm neigh}(0,\comp^2)$, smooth in
$\eps,\frac{h}{\eps^2}\in {\rm neigh}(0,\real)$, and as
$h\rightarrow 0$, has an asymptotic expansion in the space of such
functions,

$$
\widehat{P}(\xi,\eps,\frac{h}{\eps^2};h)\sim
p(\xi_1)+\eps^2\left(r_0\left(\xi,\eps,\frac{h}{\eps^2}\right)+hr_1\left(\xi,\eps,\frac{h}{\eps^2}\right)+\ldots\right).
$$
We have

$$
r_0=\langle{s}\rangle(\xi)+{\cal
O}\left(\eps+\frac{h}{\eps^2}\right),\quad r_j={\cal
O}\left(\eps+\frac{h}{\eps^2}\right),\quad j\geq 1.
$$

In the case when the subprincipal symbol of $P_{\eps=0}$ vanishes,
we have a correspon\-ding result for $\eps$ in the range $h\ll
\eps={\cal O}(h^{\delta})$, $\delta>0$, with the only difference
that "$h/\eps^2$" in the description of the eigenvalues should be
replaced by "$h/\eps$".

\end{theo}

\section{Spectral asymptotics in the torus case II}\label{section4}
\setcounter{equation}{0}

Throughout this section, it will be assumed that the subprincipal
symbol of $P_{\eps=0}$ vanishes and that $\langle{q}\rangle\equiv
0$.

Recall that the principal symbol of the original operator
$P_{\eps}$ has the form
\begeq \label{4.1}
 p_{\eps}=p+i\eps q+\eps^2 r+i\eps^3 w+{\cal O}(\eps^4),
\endeq
in a complex \neigh{} of $p^{-1}(0)\cap T^*M$. In some of the
applications that we have in mind (see section 6), the leading
perturbation $q$ is real on the real domain, and in this section
we shall assume in addition that the next order term, $r$, is real
as well. The second averaged correction $\langle{s}\rangle$,
introduced in (\ref{2.6}), is then purely real, and the results of
the preceding section do not apply. The purpose of this section is
to understand what kind of spectral results it is possible to
obtain then. As before, our starting point will be the averaging
method. To be precise, we shall now have to modify further the
choice of the IR--deformation, given by the weight $G$ in
(\ref{2.4}). Let us consider \begeq \label{4.2} G=G_0+i\eps
G_1+\eps^2 G_2+{\cal O}(\eps^3),
\endeq
where $G_2$ is to be determined, and as before, compute the
restriction of $p_{\eps}$ to the corresponding \mfld{}
$\Lambda_{\eps G}=\exp(i\eps H_G)(T^*M)$. We get
\begin{eqnarray*}
& & (p+i\eps q+\eps^2 r+i\eps^3 w)(\exp(i\eps H_G)) = p+i\eps
q+\eps^2
r+i\eps^3 w+i\eps H_G(p+i\eps q+\eps^2 r)\\
& & +\frac{(i\eps)^2}{2} H_G^2 (p+i\eps q)+\frac{(i\eps)^3}{6}H_{G_0}^3 p+ {\cal O}(\eps^4)=p+i\eps q+\eps^2
r +i\eps^3 w\\
& & +i\eps\left(H_{G_0}p+i\eps H_{G_0}q + \eps^2 H_{G_0}r+i\eps H_{G_1}p-\eps^2 H_{G_1}q+\eps^2 H_{G_2}p\right)\\
& & -\frac{\eps^2}{2}\left(H^2_{G_0}p+i\eps
  H^2_{G_0}q+i\eps(H_{G_0}H_{G_1}+H_{G_1}H_{G_0})p\right)+
\frac{(i\eps)^3}{6}H^3_{G_0}p+{\cal O}(\eps^4) \\
& = & p+\eps^2\left(-H_{G_0}q-H_{G_1}p+r-\frac{1}{2} H^2_{G_0}p \right) \\
& & +i\eps^3\left(w+H_{G_0}r-H_{G_1}q+H_{G_2}p-\frac{1}{6}
H^3_{G_0}p-\half
  H^2_{G_0}q-\half(H_{G_0}H_{G_1}+H_{G_1}H_{G_0})p\right) \\
& & +{\cal O}(\eps^4).
\end{eqnarray*}
As before, using that $H_p G_0=q$, we see that the expression in front of $\eps^2$ takes the form
$$
r+H_p G_1-\half H_{G_0}q,
$$
and the coefficient in front of $i\eps^3$ becomes
$$
w-\half H_{G_1}q-H_p G_2-\frac{1}{3} H^2_{G_0}q-\half
H_{G_0}H_{G_1}p+H_{G_0}r.
$$
From section 2 we then recall that we choose $G_1$ as an analytic
solution in a \neigh{} of $p^{-1}(0)\cap T^*M$ of the equation
$$
H_p G_1=\half\left(H_{G_0}q-\langle{H_{G_0}q\rangle}\right)-r+\langle{r}\rangle,
$$
Notice also that $G_1$ is real-valued on the real domain, since we
assumed that this is the case for $q$ and $r$.
Then we see that the
coefficient in front of $i\eps^3$ is of the form
$$
t:=w-\half H_{G_1}q -H_p
G_2-\frac{1}{12}H^2_{G_0}q-\frac{1}{4}H_{G_0}\langle{H_{G_0}q\rangle}+\half
H_{G_0}r+\half H_{G_0}\langle{r}\rangle,
$$
and it is clear that we can choose $G_2$, defined near
$p^{-1}(0)\cap T^*M$, such that after composing the leading symbol
$p_{\eps}$, which has the form (\ref{4.1}), with the complex
cano\-ni\-cal transformation $\exp(i\eps H_G)$, with $G$ of the
form (\ref{4.2}), we get a reduction of $p_{\eps}$ to \begeq
\label{4.2.1}
p+\eps^2\langle{s}\rangle+i\eps^3\langle{t}\rangle+{\cal
O}(\eps^4),
\endeq
with $\langle{s}\rangle$ and $\langle{t}\rangle$ real-valued on
the real domain, if we also assume that $w$ has this property. We
shall now give explicit expressions for these corrections,
assuming for simplicity that the period for the $H_p$--flow,
$T(E)=T$ is independent of the energy $E\in {\rm neigh}(0,\real)$.
The expression (\ref{2.6}) for $\langle{s}\rangle$ from section 2
then simplifies and we get
$$
\langle{s}\rangle=\langle{r}\rangle-\half\langle{H_{G_0}
q\rangle}=\langle{r}\rangle- \frac{1}{2T} \int_0^{T}
\{G_0,q\}\circ \exp(uH_p)\,du,
$$
and recalling the expression for $G_0$ we get \begeq \label{4.2.2}
\langle{s}\rangle=\langle{r}\rangle-\frac{1}{2T^2}\int\!\!\!\int
u\{q\circ\exp((v+u)H_p),q\circ\exp(vH_p)\}\,dv\,du,
\endeq
where the integration is performed over the rectangle $[0,T]^2$.
We also get
\begin{eqnarray*}
\langle{t}\rangle & = &
\langle{w}\rangle-\half\langle{H_{G_1}q}\rangle-
\frac{1}{12}\langle{H^2_{G_0}q}\rangle-\frac{1}{4}\langle{H_{G_0}\langle{H_{G_0}
q}\rangle}\rangle \\
& & +\half \langle{H_{G_0}r}\rangle+\half
\langle{H_{G_0}\langle{r}\rangle\rangle},
\end{eqnarray*}
where
$$
G_1=\frac{1}{2T} \int_0^{T} t \left(\left(H_{G_0} q-2r\right)\circ
\exp(tH_p)\right)\,dt\quad \wrtext{on}\quad p^{-1}(E)\cap T^*M.
$$
We shall now simplify the expression for $\langle{t}\rangle$. First of all, notice that
$$
\langle{H_{G_0}\langle{H_{G_0} q}\rangle}\rangle=\frac{1}{T}
\int_0^T \{G_0\circ \exp(tH_p),\langle{H_{G_0}q}\rangle\}\,dt=\{\langle{G_0}\rangle,\langle{H_{G_0}q\rangle}\}=0,
$$
since we have already observed that the flow average of $G_0$
vanishes. Similarly,
$$
\langle{H_{G_0}\langle{r}\rangle\rangle}=0.
$$
If we introduce \begeq \label{4.2.3.1} f=-\frac{1}{12T^2}
\int\!\!\!\int uv \{q\circ \exp(vH_p),\{q\circ \exp(u
H_p),q\}\}\,du\,dv,
\endeq
\begeq \label{4.2.3.2} g=-\frac{1}{4T^2} \int\!\!\!\int uv
\{\{q\circ \exp((u+v)H_p), q\circ \exp(uH_p)\},q\}\,du\,dv,
\endeq
and \begeq \label{4.2.3.3}
 k=\frac{1}{2T}\int_0^T u\left(\{q\circ
\exp(uH_p),r\}+\{r\circ\exp(uH_p),q\}\right)\,du,
\endeq
then \begeq \label{4.2.4}
\langle{t}\rangle=\langle{f+g+w+k}\rangle.
\endeq

\vskip 4mm In what follows, we shall consider, as we may, the
operator $P_{\eps}$, defined microlocally near $p^{-1}(0)\cap
T^*M$, with the principal symbol of the form
$$
p+\eps^2 \langle{s}\rangle+i\eps^3 \langle{t}\rangle+{\cal O}(\eps^4),
$$
with $\langle{s}\rangle$ and $\langle{t}\rangle$ real-valued, and
with the subprincipal symbol ${\cal O}(\eps)$. When deriving
spectral asymptotics for $P_{\eps}$, we shall distinguish the
following cases, depending on the size of $\eps$:
\begin{enumerate}
\item $h \ll \eps \ll h^{1/2}$. \item $h^{1/2} \ll \eps={\cal
O}(h^{\delta}),\quad \delta>0$. \item $\eps \sim h^{1/2}$.
\end{enumerate}

\medskip \noindent 1. {\it The case when $h\ll \eps \ll h^{1/2}$}.

In this case we know, thanks to Proposition 2.1, that the spectrum
of $P_{\eps}$ near $0$ has a cluster structure.

Let us fix the energy level $(0,0)$ and introduce the set
$$
\Lambda_{0,0}: p=0,\, \langle{s}\rangle=0,
$$
which is invariant under the $H_p$-flow. We assume that $dp$ and
$d\langle{s}\rangle$ are linearly independent at every point of
$\Lambda_{0,0}$, and then, assuming that $\Lambda_{0,0}$ is
connected, we see that this set is a Lagrangian torus. We shall
also assume that $T(0)$ is the minimal period of every closed
$H_p$--trajectory in $\Lambda_{0,0}$. In a \neigh{} of this torus,
$p$ and $\langle{s}\rangle$ form a completely integrable system,
and we then pass to the action-angle variables by means of a real
analytic canonical transformation
$$
\kappa: {\rm neigh}\left(\xi=0, T^*{\bf T}^2\right)\rightarrow
{\rm neigh}\left(\Lambda_{0,0},T^*M\right),
$$
such that $p\circ\kappa=p(\xi_1)$,
$\langle{s}\rangle\circ\kappa=\langle{s}\rangle(\xi)$. Notice that
when expressed in terms of the action coordinate $\xi_1$, $p$
becomes $p(\xi_1)=f(\xi_1)$, where the function $f$ has been
defined after (\ref{2.6.5}). Implementing $\kappa$ by means of a
microlocally unitary \fourior{}, we get a new operator, still
denoted by $P_{\eps}$, microlocally defined near the zero section
in $T^*{\bf T}^2$, which has the leading symbol
$$
p(\xi_1)+\eps^2 \langle{s}\rangle(\xi)+i\eps^3 \langle{t}\rangle(x_2,\xi)+{\cal O}(\eps^4).
$$
The assumption about the linear independence of the differentials
implies that
$$
\partial_{\xi_1}p(0)\neq 0, \quad \partial_{\xi_2}\langle{s}\rangle(0)\neq
0,
$$
and repeating the arguments of section 3, we then find an
$h$-\pseudor{} $A=\sum_{\nu=0}^{\infty} h^{\nu}
a_{\nu}(x,\xi,\eps)$, with $a_0={\cal O}(\eps^4)$, such that the
full symbol of the operator
$$
e^{\frac{i}{h} A}P_{\eps} e^{-\frac{i}{h} A}=e^{\frac{i}{h}{\rm ad}
  A}P_{\eps}=\sum_{k=0}^{\infty}
\left(\frac{i}{h}{\rm ad}\,A\right)^k P_{\eps}
$$
is independent of $x_1$. This is done by solving transport
equations in $x_1$, obtained by looking for each
$a_{\nu}(x,\xi,\eps)$ as a formal power series in $\eps$. This
reduces the following discussion to the operator
$$
P_{\eps}=\sum_{j=0}^{\infty} h^j p_j(x_2,\xi,\eps),
$$
with
\begeq
\label{final}
p_0=p(\xi_1)+\eps^2 \langle{s}\rangle(\xi)+i\eps^3\langle{t}\rangle(x_2,\xi)+{\cal O}(\eps^4)
\endeq
also independent of $x_1$, and $p_1={\cal O}(\eps)$. It is then clear that under the assumption that
$$
\frac{h}{\eps}\leq h^{\delta_1},\quad \delta_1>0,
$$
by repeating the arguments of section 3, we can perform further
conjugations of $P_{\eps}$ by means of Fourier integral operators,
and make its full symbol independent of $x_2$, modulo ${\cal
O}(h^{\infty})$. To be precise, we find that there exist
$$
B_0=b_0(x_2,hD_x,\eps,\frac{h}{\eps};h)\quad b_0={\cal
O}\left(\eps+\frac{h}{\eps}\right),
$$
and
$$
B_1=\sum_{\nu=1}^{\infty} h^{\nu}
b_{\nu}(x_2,hD_x,\eps,\frac{h}{\eps}),\quad b_{\nu}={\cal
O}\left(\eps+\frac{h}{\eps}\right),
$$
such that
$$
\widehat{P}=e^{\frac{i}{h}{\rm ad}_{B_1}}e^{\frac{i}{h}{\rm
ad}_{B_0}}{P}_{\eps}
$$
has the full symbol independent of $x$, modulo ${\cal
O}(h^{\infty})$. In the case when we only have $h/\eps \ll 1$, as
before, the operator $e^{iB_0/h}$ has to be replaced by a more
general Fourier integral operator, constructed as in section 3 and
as in~\cite{HiSj}. We then have the associated quasi-eigenvalues
given by
$$
\widehat{P}\left(h(k-\theta),\eps,\frac{h}{\eps};h\right),\quad k\in
\z^2,\quad \theta=\frac{S}{2\pi h}+\frac{\alpha}{4}.
$$
Analyzing this microlocal reduction to a normal form near
$\Lambda_{0,0}$, similarly to what we did in section 3 and in
section 6 of~\cite{HiSj}, we get the following result.

\begin{prop}
There exists an IR-manifold $\Lambda\subset T^*\widetilde{M}$,
which is $(\eps+h/\eps)$-close to $T^*M$ and agrees with $T^*M$
away from $p^{-1}(0)\cap T^*M$, and a smooth Lagrangian torus
$\widetilde{\Lambda}_{0,0}\subset \Lambda$, such that for $\rho
\in \Lambda$ in a \neigh{} of $\widetilde{\Lambda}_{0,0}$, we have
$$
\rho=\exp(i\eps H_G)\circ\kappa\circ \widetilde{\kappa}(x,\xi),\quad
(x,\xi)\in T^*{\bf T}^2,
$$
and $\widetilde{\Lambda}_{0,0}$ corresponds precisely to the zero
section $\xi=0$. Here $\kappa$ is the canonical transformation
given by the action-angle variables, and $\widetilde{\kappa}$ is
constructed as before, using the generating function {\rm
(\ref{3.6})}, with $h/\eps^2$ replaced there by $h/\eps$. When
$\rho\in \Lambda$ is away from a small \neigh{} of
$\widetilde{\Lambda}_{0,0}$ in $\Lambda$ and $\abs{\Re
  P_{\eps}(\rho;h)}\leq 1/C$ for a large enough $C$, it is true that
$$
\abs{\langle{s}\rangle(\rho)}\geq \frac{1}{{\cal O}(1)}.
$$
Furthermore, there exists an elliptic uniformly bounded \fourior{} $U:
H(\Lambda)\rightarrow L^2_{\theta}({\bf T}^2)$ such that, microlocally
near $\widetilde{\Lambda}_{0,0}$, $UP_{\eps}=\widehat{P} U$. Here
$$
\widehat{P}=\widehat{P}\left(hD_x,\eps,\frac{h}{\eps};h\right): L^2_{\theta}({\bf
  T}^2)\rightarrow L^2_{\theta}({\bf T}^2)
$$
has the full symbol holomorphic in $\xi\sim 0$ in $\comp^2$, and
depending smoothly on $\eps$, $h/\eps\in {\rm neigh}(0,\real)$,
$$
\widehat{P}(\xi,\eps,\frac{h}{\eps};h)=p(\xi_1)+\eps^2
\left(r_0(\xi,\eps, \frac{h}{\eps})+hr_1(\xi,\eps,\frac{h}{\eps})+\ldots\right),
$$
with
$$
r_0=\langle{s}\rangle(\xi)+{\cal
  O}\left(\eps+\frac{h}{\eps}\right),\quad r_j={\cal O}\left(\eps+\frac{h}{\eps}\right),\quad j\geq
1.
$$
\end{prop}

Before setting up a suitable globally well-posed Grushin problem in
the Hilbert space $H(\Lambda)$ in order to identify the
spectrum precisely, following~\cite{HiSj}, we shall first show
that when $z\in {\rm neigh}(0,\comp)$ is such that \begeq
\label{0.31} \abs{\Re
z-f\left(h\left(k_1-\frac{\alpha_1}{4}\right)-\frac{S_1}{2\pi}\right)}\leq
\frac{\eps^2}{C}, \quad C\gg 1,
\endeq
for some $k_1\in \z$, and $z$ avoids the union of the pairwise
disjoint open discs $D_{k_2}(h)$ of radii $\eps^2 h/{\cal O}(1)$,
that are centered at the quasi--eigenvalues
$$
\widehat{P}
\left(h(k_1-\frac{\alpha_1}{4})-\frac{S_1}{2\pi},h(k_2-\frac{\alpha_2}{4})-\frac{S_2}{2\pi},\eps,\frac{h}{\eps};h\right),
$$
for $k_2\in \z$, then the operator
$$
P_{\eps}-z: H(\Lambda) \rightarrow H(\Lambda)
$$
is bijective. In doing so, we shall essentially follow an argument
from section 6 of~\cite{HiSj}, with a slight improvement. Let
$\chi\in C^{\infty}_0(\Lambda)$ be supported in a small flow
invariant \neigh{} of $\widetilde{\Lambda}_{0,0}$ where $P_{\eps}$
is intertwined with $\widehat{P}$, and $\chi=1$ near
$\widetilde{\Lambda}_{0,0}$. We also assume, as we may, that on
the operator lever,
\begeq \label{0.31.1} [P_{\eps}, \chi]={\cal
O}(h^{\infty}): H(\Lambda)\rightarrow H(\Lambda).
\endeq
Indeed, to achieve the property (\ref{0.31.1}) it suffices to
choose a function $\chi_0\in C^{\infty}_0(T^*\T^2)$ with
$\supp\,\chi_0$ close to $\xi=0$ and $\chi_0=1$ in a \neigh{} of
this set, such that $\chi_0=\chi_0(\xi)$ depends on $\xi$ only.
Conjugating $\chi_0$ by means of the microlocal inverse $V$ of $U$
and modifying it by an ${\cal O}(h^{\infty})$-factor outside a
\neigh{} of $\widetilde{\Lambda}_{0,0}$, we obtain $\chi$ with the
desired properties.

Let us introduce a partition of unity on $\Lambda$,
$$
1=\chi+\psi_{1,+}+\psi_{1,-}+\psi_{2,+}+\psi_{2,-}.
$$
Here the functions $\psi_{1,\pm}\in C^{\infty}_0(\Lambda)$ are
supported in flow invariant regions $\Omega_{\pm}$, such that
$\pm \langle{s}\rangle> 1/{\cal O}(1)$ in $\Omega_{\pm}$,
respectively. Moreover, we can arrange so that $\psi_{1,\pm}$ are in involution with
$p$, the principal symbol of $P_{\eps=0}$ on $H(\Lambda)$. Finally, $\psi_{2,\pm}$ are such that $\pm \Re P_{\eps}>1/{\cal
O}(1)$ in the support of $\psi_{2,\pm}$. We shall prove that
\begeq
\label{0.32}
\norm{(1-\chi)u}\leq {\cal
O}\left(\frac{1}{\eps^2}\right)\norm{v}+{\cal O}(h^{\infty})\norm{u}.
\endeq
Here the norms are taken in $H(\Lambda)$, $(P_{\eps}-z)u=v$, and $z\in
\comp$ satisfies (\ref{0.31}) and avoids the discs
$D_{k_2}(h)$, $k_2\in \z$. When establishing (\ref{0.32}), we shall first prove
this bound with $\psi_{1,+}$ in place of $1-\chi$.

When $N\in \nat$, let $\psi_0\prec \psi_1\prec\ldots \prec
\psi_{N}$, $\psi_{0}=\psi_{1,+}$, be smooth cutoff functions
supported in $\Omega_+$, which Poisson commute with $p$. (The
standard notation $f\prec g$ means that $g=1$ in a \neigh{} of
$\supp\,f$.) We have
$$
(P_{\eps}-z)\psi_j u=\psi_j v+[P_{\eps},\psi_j]u,\quad 0\leq j\leq N.
$$
Now the operator $P_{\eps}-z$ is invertible, microlocally in
$\Omega_+$, with the norm of the microlocal inverse being ${\cal
O}(\eps^{-2})$. It follows that
$$
\norm{\psi_{j}u}\leq {\cal
O}\left(\frac{1}{\eps^2}\right)\left(\norm{v}+\norm{[P_{\eps},\psi_{j}]u}\right)+{\cal
O}(h^{\infty})\norm{u}.
$$
Notice that here $[P_{\eps},\psi_j]={\cal O}(h^3)+{\cal O}(\eps^2
h)={\cal O}(\eps^2 h)$, since the subprincipal symbol of $P_{\eps=0}$
vanishes and $h\leq \eps$. Now
$\psi_{j}(1-\psi_{j+1})={\cal O}(h^{\infty})$ in the
operator sense, and hence we obtain
$$
\norm{[P_{\eps},\psi_{j}]u}\leq {\cal O}(\eps^2 h)\norm{\psi_{j+1}u}+{\cal O}(h^{\infty})\norm{u},
$$
and therefore,
$$
\norm{\psi_{j}u}\leq {\cal O}\left(\frac{1}{\eps^2}\right)
\norm{v}+{\cal O}(h)\norm{\psi_{j+1}u}+{\cal O}(h^{\infty})\norm{u}.
$$
Combining these estimates for $j=0,\ldots N$, we get
$$
\norm{\psi_{1,+}u}\leq {\cal O}\left(\frac{1}{\eps^2}\right)
\norm{v}+{\cal O}_N(1) h^{N}\norm{\psi_N u}+{\cal O}(h^{\infty})\norm{u}.
$$
The same estimate is then obtained for $\psi_{1,-}u$, which is
concentrated in the region where $\langle{s}\rangle<-1/{\cal O}(1)$, and
a fortiori such estimates also hold in regions where $\pm \Re P_{\eps}\sim 1$. The
estimate (\ref{0.32}) follows.

Relying upon (\ref{0.32}), we  shall complete the proof of the fact
that the spectrum of $P_{\eps}$ in the set (\ref{0.31}) is contained
in the union of the discs $D_{k_2}(h)$. Let us write
$$
(P_{\eps}-z)\chi u=\chi v+[P_{\eps},\chi]u,
$$
and from (\ref{0.31.1}) we recall that the norm of the commutator term
does not exceed
$$
{\cal O}(h^{\infty})\norm{u}.
$$
Applying the \fourior{} $U$ of Proposition 4.1, we get, modulo
${\cal O}(h^{\infty})$-errors,
$$
\left(\widehat{P}-z\right)U\chi u=U\left(\chi v+[P_{\eps},\chi]u\right).
$$
Now an expansion in Fourier series shows that the operator $\widehat{P}-z$ is invertible, microlocally near
$\xi=0$, with a microlocal inverse of the norm ${\cal
O}(\eps^{-2}h^{-1})$, provided that $z$ in the set (\ref{0.31}) avoids the discs $D_{k_2}(h)$. We get
$$
\norm{\chi u}\leq {\cal O}\left(\frac{1}{\eps^2 h}\right)\norm{v}+{\cal O}(h^{\infty})\norm{u},
$$
and combining this estimate together with (\ref{0.32}) we infer
that the operator $P_{\eps}-z~:H(\Lambda)\rightarrow H(\Lambda)$
is injective, hence bijective, since it is a Fredholm operator of
index zero by general arguments. Now $H(\Lambda)$ agrees with
$L^2(M)$ as a space, and we conclude that $z$ is not in the
spectrum of $P_{\eps}$.

When $z$ varies in the disc $D_{k_2}(h)$ contained in the set
(\ref{0.31}), for some $k_2\in \z$, we shall recall briefly the
setup of the global Grushin problem for the operator $P_{\eps}-z$.
Introduce the operators $R_+: H(\Lambda)\rightarrow \comp$ and
$R_-: \comp\rightarrow H(\Lambda)$ by
$$
R_+u=(U\chi u|e_k),\quad R_-u_-=u_-Ve_k,\quad k=(k_1,k_2).
$$
Here
$$
e_k=\frac{1}{2\pi}
\exp\left(\frac{i}{h}(h(k-\frac{\alpha}{4})-\frac{S}{2\pi})x\right),
$$
the scalar product in the definition of $R_+$ is taken in
$L^2_{\theta}(\T^2)$, and $V$ is the microlocal inverse of $U$. It
is then clear that the arguments of section 6 of~\cite{HiSj} apply
as they stand and show that for every $(v,v_+)\in H(\Lambda)\times
\comp$, the Grushin problem
$$
(P_{\eps}-z)u+R_-u_-=v,\quad R_+u=v_+
$$
has a unique solution $(u,u_-)\in H(\Lambda)\times \comp$ given by
$$
u=Ev+E_+v_+,\quad u_-=E_-v+E_{-+}v_+.
$$
Here the zeros of $E_{-+}(z)$ agree with the eigenvalues of $P_{\eps}$ in
$D_{k_2}(h)$, and repeating the arguments of~\cite{HiSj} we find that $E_{-+}(z)=z-\widehat{P}(h(k-\theta_1),
h(k_2-\theta_2),\eps,\frac{h}{\eps};h)$ modulo an error term which is
${\cal O}(h^{\infty})$. We get the following result.

\begin{theo}
Let $P_{\eps}$ be an operator as in the introduction, so that {\rm
(\ref{1.12})}, {\rm (\ref{4.1})} hold. Assume that $q$, $r$, $w$
are real-valued and that the trajectory average
$\langle{q}\rangle$ vanishes. Define $s$ by {\rm (\ref{2.6})} and
recall that $p_{\eps}$ can be reduced by averaging to {\rm
(\ref{4.2.1})}, where $\langle{t}\rangle$ is given by {\rm
(\ref{4.2.4})}, {\rm (\ref{4.2.3.1})}--{\rm (\ref{4.2.3.3})}, in
the case when $T(E)=T$ is independent of $E$. Assume that the
subprincipal symbol of $P_{\eps=0}$ vanishes. Let $F_0$ be a
regular value of $\langle{s}\rangle$, considered as a function on
$p^{-1}(0)\cap T^*M$. Assume that the Lagrangian manifold
$$
\Lambda_{0,F_0}: p=0,\,\,\langle{s}\rangle=F_0
$$
is connected and that $T(0)$ is a minimal period of every closed
$H_p$-trajectory in $\Lambda_{0,F_0}$. When $\gamma_1$ and
$\gamma_2$ are the fundamental cycles in $\Lambda_{0,F_0}$ with
$\gamma_1$ being given by a closed $H_p$--trajectory, we write
$S=(S_1,S_2)$ and $\alpha=(\alpha_1,\alpha_2)$ for the actions and
Maslov indices of the cycles, respectively. Assume that
$$
h\ll \eps \ll h^{1/2}.
$$
Let $C>0$ be large enough. Then for each $k_1\in \z$, the eigenvalues of
$P_{\eps}$ in the set
$$
\abs{\Re
z-f\left(h(k_1-\frac{\alpha_1}{4})-\frac{S_1}{2\pi}\right)-\eps^2 F_0}<
\frac{\eps^2}{C}
$$
are given by
\begeq
\label{0.33}
\widehat{P}\left(h(k_1-\frac{\alpha_1}{4})-\frac{S_1}{2\pi},h(k_2-\frac{\alpha_2}{4})-\frac{S_2}{2\pi},\eps,\frac{h}{\eps};h\right)
+{\cal O}(h^{\infty}),\quad k_2\in \z.
\endeq
Here $\widehat{P}(\xi,\eps,h/\eps;h)$ is holomorphic in $\xi \in {\rm
neigh}(0,\comp^2)$, smooth in $\eps,\frac{h}{\eps}\in {\rm
neigh}(0,\real)$, and has a complete asymptotic expansion in the
space of such functions, as $h\rightarrow 0$,
$$
\widehat{P}\left(\xi,\eps,\frac{h}{\eps};h\right)\sim
f(\xi_1)+\eps^2\left(r_0\left(\xi,\eps,\frac{h}{\eps}\right)+h
r_1\left(\xi,\eps,\frac{h}{\eps}\right)+\ldots\right).
$$
We have
$$
r_0(\xi)=\langle{s}\rangle(\xi)+{\cal
O}\left(\eps+\frac{h}{\eps}\right),\quad r_j={\cal O}(\eps+\frac{h}{\eps}),\quad j\geq
1.
$$
\end{theo}

\Remark. It is understood that in (\ref{0.33}) we only consider those
values of $k_1\in \z$ for which the first argument ($\xi_1$) of
$\widehat{P}$ is small enough.

\Remark. The arguments presented in this case can also be applied
to a class of {\it selfadjoint} perturbations of $P_{\eps=0}$,
with the leading symbol of the form $p+\eps q+{\cal O}(\eps^2 m)$,
and with the subprincipal symbol ${\cal O}(\eps m)$. (In the
compact case, the weight $m$ should be replaced by
$\langle{\xi}\rangle^m$.) Here $q$ is real on $T^*M$ and we assume
that the flow average $\langle{q}\rangle$ does not vanish
identically. If $\eps \ll h$ then we know from Proposition 2.1
that the spectrum of $P_{\eps}$ in ${\rm neigh}(0,\comp)$ has a
cluster structure. Let us restrict $\eps$ further and assume that
$h^2 \ll \eps$. If $F_0$ is a regular value of $\langle{q}\rangle$
along $p^{-1}(0)$, then it follows from the arguments of this
subsection that we may describe precisely the asymptotics of the
individual eigenvalues of $P_{\eps}$ inside the sub-clusters given
by
$$
\abs{\Re
z-f\left(h\left(k_1-\frac{\alpha}{4}\right)-\frac{S_1}{2\pi}\right)-\eps
F_0}<\frac{\eps}{C},\quad C\gg 1,\quad k_1\in \z.
$$
Notice, in particular, that if $M$ is a smooth Zoll surface, and
$P=-\Delta+\lambda V$ is a Schr\"odinger operator on $M$, then the
operator $h^2P$ is of the form above, provided that the coupling
constant $\lambda$ in front of the potential is sufficiently
large, but fixed. Also, notice that in this case, the analyticity
assumptions are not needed, as we may work with Fourier integral
operators with real phase and hence stay within the framework of
the usual $L^2$--spaces. To the best of our knowledge, such
results in the selfadjoint case do not seem to have been
formulated explicitly in the literature---see, however,~\cite{UZ}
for asymptotics of pair correlation functions on Zoll surfaces.

\medskip \noindent 2. {\it The case when $h^{1/2}\ll \eps={\cal
O}(h^{\delta}),\quad \delta>0$.}

Let us recall that by means of an averaging procedure, we have
reduced the original ope\-ra\-tor $P_{\eps}$ to an operator with
the principal symbol \begeq \label{0.35} p+\eps^2
\langle{s}\rangle+i\eps^3 \langle{t}\rangle+{\cal O}(\eps^4),
\endeq
where $\langle{s}\rangle$ and $\langle{t}\rangle$ are real on the real
domain. Moreover, the subprincipal symbol of the averaged operator is
${\cal O}(\eps)$.

In the first case, we saw that due to the cluster structure of the
spectrum, the imaginary part of the principal symbol played no
major role in the arguments. The situation will be quite different
now. When $T>0$, we introduce the double average, \begeq
\label{0.36} \langle{\langle{t}\rangle}\rangle_T =\frac{1}{T}
\int_0^T \langle{t}\rangle\circ \exp(uH_{\langle{s}\rangle})\,du,
\endeq
defined in a \neigh{} of $p^{-1}(0)\cap T^*M$. Let $G=G_T$ be a
real-valued analytic function defined near $p^{-1}(0)\cap T^*M$,
such that \begeq \label{0.37} H_{\langle{s}\rangle}G =
\langle{t}\rangle-\langle{\langle{t}\rangle}\rangle_T.
\endeq
The equation (\ref{0.37}) is a convolution equation along the
$H_{\langle{s}\rangle}$-trajectories, and we solve it by putting
$$
G = \int k(u/T) \langle{t}\rangle\circ
\exp(uH_{\langle{s}\rangle})\,du,
$$
where $k$ is a piecewise linear function, vanishing outside a
bounded interval, and solving $k' =-\delta_0+1_{[0,1]}$. Notice
that $H_pG=0$. After composing the principal symbol (\ref{0.35})
with the holomorphic \canform{} $\exp(i\eps H_G)$, and, on the
operator level, after a conjugation by means of the corresponding
microlocally unitary \fourior{} which also has the improved Egorov
property, we may assume that our operator is microlocally defined
near $p^{-1}(0)\cap T^*M$ and has the principal symbol \begeq
\label{0.38} p+\eps^2
\langle{s}\rangle+i\eps^3\langle{\langle{t}\rangle}\rangle_T+
{\cal O}_T(\eps^4).
\endeq
We also know that the subprincipal symbol of the new operator is
${\cal O}_T(\eps)$, and we can be more precise and remark that it is
in fact ${\cal O}(\eps)+{\cal O}_T(\eps^2)$.

We keep the basic assumption of case 1 that $dp$ and
$d\langle{s}\rangle$ are linearly independent along the set
$\Lambda_{0,0}: p=0, \langle{s}\rangle=0$, which is then a
flow-invariant Lagrangian torus, provided that it is connected.
Assume also that $T(0)$ is the minimal period for the $H_p$-flow
in $\Lambda_{0,0}$. Consider now the restriction of
$\langle{\langle{t}\rangle}\rangle_T$ to $\Lambda_{0,0}$. Then the
limit \begeq \label{0.39}
\langle\langle{t}\rangle\rangle_{\infty}~:=\lim_{T\rightarrow
\infty} \langle\langle{t}\rangle\rangle_T(\rho),\quad \rho \in
\Lambda_{0,0},
\endeq
exists, $\langle{\langle{t}\rangle}\rangle_{\infty}\in \real$, and the
  convergence is uniform on $\Lambda_{0,0}$. We also remark that
  $\langle{\langle{t}\rangle}\rangle_{\infty}$ is precisely the mean value of
  $\langle{t}\rangle$ over $\Lambda_{0,0}$ and that
\begeq
\label{0.40}
\langle\langle{t}\rangle\rangle_T(\rho)-\langle\langle{t}\rangle\rangle_{\infty}={\cal
  O}\left(\frac{1}{T}\right),
\endeq
uniformly on $\Lambda_{0,0}$. Let us also recall the Lagrangian
foliation given by the flow-invariant tori
$$
\Lambda_{E,F}: p=E,\,\,\langle{s}\rangle=F,
$$
for $(E,F)\in {\rm neigh}(0,\real^2)$. We shall
assume that for $0<a\ll 1$, the set
$$
p^{-1}(0)\backslash \bigcup_{\abs{F}<a}\Lambda_{0,F}
$$
splits into two disjoint flow-invariant components, which we
denote by $\Lambda_a$ and $\Lambda_{-a}$. We then introduce the
following global assumption:
\begin{eqnarray}
\label{H} & & \hbox{For any}\,\,b\in (0,a)\,\,\hbox{there exists}
\\ \nonumber & & T(b)>0\,\,\hbox{and}\,\,C(b)>0\,\,\hbox{such
that}\\ \nonumber & &
\inf_{\Lambda_{[b,a]}}\biggl(\langle{\langle{t}\rangle}\rangle_T(\rho)-\langle\langle{t}\rangle\rangle_{\infty}\biggr)
\geq \frac{1}{C(b)},\quad T\geq T(b),\\ \nonumber & &
\sup_{\Lambda_{[-a,-b]}}\biggl(\langle{\langle{t}\rangle}\rangle_T(\rho)-\langle\langle{t}\rangle\rangle_{\infty}\biggr)
\leq -\frac{1}{C(b)},\quad T\geq T(b).
\end{eqnarray}
Here we have put
$$
\Lambda_{[b,a]}= \left (\bigcup_{b\leq
F<a}\Lambda_{0,F}\right)\bigcup \Lambda_{a},
$$
and $\Lambda_{[-a,-b]}$ is defined similarly.

Given an operator with the principal symbol (\ref{0.38}), considered
microlocally near $\Lambda_{0,0}$, we pass to the standard torus ${\bf
T}^2$ and eliminate the $x_1$-variable, exactly as before. We are then
reduced to the operator
$\widetilde{P}_{\eps}=\widetilde{P}_{\eps}(T)$, acting on
$L^2_{\theta}(\T^2)$, which is defined microlocally
near $\xi=0$ in $T^*{\bf T}^2$, and which has the form
\begeq
\label{0.43}
\widetilde{P}_{\eps}=\sum_{j=0}^{\infty}
h^{j}\widetilde{p}_j(x_2,\xi,\eps,T),
\endeq
with
$$
\widetilde{p}_0(x_2,\xi,\eps,T)=p(\xi_1)+\eps^2 \langle{s}\rangle(\xi)+i\eps^3
\langle\langle{t}\rangle\rangle_T(x_2,\xi)+{\cal O}_T(\eps^4),
$$
and
$$
\widetilde{p}_1(x_2,\xi,\eps,T)=\eps \widetilde{q}_1(x_2,\xi,\eps,T)=
\eps\widetilde{q}_1(x_2,\xi,0)+{\cal O}_T(\eps^2).
$$
Let us write \begeq \label{0.44} \widetilde{P}_{\eps} =
p(\xi_1)+\eps^2
\left(r_0(x_2,\xi,\eps,\frac{h}{\eps},T)+hr_1(x_2,\xi,\eps,\frac{h}{\eps},T)+\ldots\right),
\endeq
where
$$
r_0=\langle{s}\rangle(\xi)+i\eps\langle\langle{t}\rangle\rangle_T(x_2,\xi)+{\cal
  O}_T(\eps^2)+
\frac{h}{\eps}\widetilde{q}_1(x_2,\xi,0)+\frac{h^2}{\eps^2}p_2(x_2,\xi,\eps,T),
$$
and
$$
r_1={\cal O}_T(1)+\frac{h^2}{\eps^2}\widetilde{p_3},
$$
$$
r_2=\frac{h^2}{\eps^2}\widetilde{p_4},\ldots,
$$
so that $r_j={\cal O}_T(1)$, $j\geq 1$.

We shall reexamine the construction of the operator
$B_0=b_0(x_2,hD_x,\eps,\frac{h}{\eps},T)$ such that
$$
r_0\circ \exp(H_{b_0})=\sum_{k=0}^{\infty} \frac{H_{b_0}^k r_0}{k!}
$$
is independent of $x_2$, modulo ${\cal O}(h^{\infty})$. As before,
we are looking for $b_0$ in terms of a formal power series in
$h/\eps={\cal O}(h^{1/2})$ and $\eps={\cal O}(h^{\delta})$,
$\delta>0$. Writing
$$
b_0=b_{0,1,0}\eps+b_{0,0,1}\frac{h}{\eps}+{\cal
  O}\left((\eps,\frac{h}{\eps})^2\right),
$$
we find that $b_{0,0,1}$ should solve the transport equation
$$
\partial_{\xi_2}\langle{s}\rangle \partial_{x_2}
b_{0,0,1}=\widetilde{q}_1(x_2,\xi,0)-\langle{\widetilde{q_1}(\cdot,\xi,0)}\rangle,
$$
with the average of $\widetilde{q}_1$ in the right hand side standing for the average with
respect to the $x_2$ variable. Therefore, $b_{0,0,1}={\cal O}(1)$. As
for the term $b_{0,1,0}$, we see that it satisfies
$$
\partial_{\xi_2}\langle{s}\rangle \partial_{x_2}b_{0,1,0}=i\left(\langle\langle{t}\rangle\rangle_T(x_2,\xi)-
\langle\langle\langle{t}\rangle\rangle_T\rangle\right),
$$
where again the last average in the right hand side is with respect to
$x_2$. Since the right hand side here is ${\cal O}(1/T)$, uniformly in $x_2$ and
$\xi\in {\rm neigh}(0,\real^2)$, we get
$$
b_{0,1,0}(x_2,\xi)={\cal O}\left(\frac{1}{T}\right).
$$
In a similar way, we construct the lower order terms, and obtain
$$
b_0={\cal O}\left(\frac{1}{T}\right)\eps+{\cal O}\left(\frac{h}{\eps}\right)+{\cal
  O}_T\left((\eps,\frac{h}{\eps})^2\right).
$$
It follows that there exists
$$
B_0=b_0(x_2,hD_x,\eps,\frac{h}{\eps},T),\quad b_0={\cal
O}\left(\frac{\eps}{T}\right)+{\cal
O}\left(\frac{h}{\eps}\right)+{\cal O}_T\left((\eps,\frac{h}{\eps})^2\right),
$$
and
$$
B_1=\sum_{\nu=1}^{\infty} h^{\nu}
b_{\nu}(x_2,hD_x,\eps,\frac{h}{\eps},T),\quad b_{\nu}={\cal
O}_T(1),
$$
such that
$$
\widehat{P}=e^{\frac{i}{h}{\rm ad}_{B_1}}e^{\frac{i}{h}{\rm ad}_{B_0}}\widetilde{P}_{\eps}
$$
has the full symbol independent of $x$,
$$
\widehat{P}=p(\xi_1)+\eps^2\left(r_0\left(\xi,\eps,\frac{h}{\eps},T\right)+h
r_1\left(\xi,\eps,\frac{h}{\eps},T\right)+\ldots\right),
$$
with
$$
r_0=\langle{s}\rangle(\xi)+i\eps\langle\langle\langle{t}\rangle\rangle_T\rangle(\xi)+{\cal
O}\left(\frac{h}{\eps}\right)+{\cal O}_T\left(\eps^2+\left(\frac{h}{\eps}\right)^2\right),
$$
$$
r_j={\cal O}_T(1),\quad j\geq 1.
$$
Repeating the arguments of section 6 of~\cite{HiSj}, we now introduce a globally defined IR-manifold
$\widetilde{\Lambda}\subset T^*\widetilde{\bf T}^2$, with
$\widetilde{\bf T}^2$ standing for the standard complexification of
${\bf T}^2$, such that $\widetilde{\Lambda}=\exp(H_{b_0})({\rm
  neigh}\,(\xi=0,T^*{\bf T}^2))$ near $\xi=0$ in $T^*\widetilde{\bf
  T}^2$, $\widetilde{\Lambda}=T^*{\bf T}^2$ further away from
$\xi=0$, and such that along $\widetilde{\Lambda}$,
\begeq
\label{0.45}
\Im \xi_1=0,\quad \Im \xi_2={\cal
  O}\left(\frac{\eps}{T}\right)+{\cal O}\left(\frac{h}{\eps}\right)+
{\cal O}_T\left((\eps,\frac{h}{\eps})^2\right),
\endeq
$$
\Im x={\cal O}\left(\frac{\eps}{T}\right)+{\cal O}\left(\frac{h}{\eps}\right)+
{\cal O}_T\left((\eps,\frac{h}{\eps})^2\right).
$$
Then the action of $\widetilde{P}_{\eps}$ on the space $H(\widetilde{\Lambda})$,
associated to $\widetilde{\Lambda}$ by means of the FBI-Bargmann
transformation on $\T^2$, is microlocally near $\xi=0$, unitarily equivalent to the conjugated
operator which acts on $L^2({\bf T}^2)$ and whose Weyl symbol is independent of $x$ and has the form
$$
p(\xi_1)+\eps^2\left(r_0(\xi,\eps,\frac{h}{\eps},T)+h
  r_1(\xi,\eps,\frac{h}{\eps},T)+\ldots\right),
$$
with the same estimates on $r_j$ as above.

We consider the imaginary part of $\widetilde{P}_{\eps}$ along
$\widetilde{\Lambda}$,
$$
\Im\widetilde{P}_{\eps}=\eps^3 \langle\langle{t}\rangle\rangle_T(\Re x_2,\Re
\xi)+{\cal O}\left(\frac{\eps^3}{T}\right)+{\cal O}(\eps h)+{\cal O}_T(\eps^4).
$$
Choosing now first $T$ sufficiently large but fixed, and then $h$
small enough, and using the fact that $h^{1/2}\ll \eps$ together
with the global assumption (\ref{H}), we see that away from a
small \neigh{} of the torus, the absolute value of $\Im
\widetilde{P}_{\eps}-\eps^3
\langle\langle{t}\rangle\rangle_{\infty}$ is bounded from below by
$\eps^3/{{\cal O}(1)}$, provided that $\abs{\Re
\widetilde{P}_{\eps}}\leq 1/{\cal O}(1)$.

Corresponding to the manifold $\widetilde{\Lambda}$ on the torus
side, we get a globally defined IR-manifold $\Lambda\subset
T^*\widetilde{M}$, which is an $(\eps+h/\eps)$-perturbation of
$T^*M$, and near $\Lambda_{0,0}$ it is obtained by replacing
$\kappa(T^*\T^2)$ there by $\kappa(\widetilde{\Lambda})$. Here
$\kappa$ is the action-angle \canform, and for simplicity here we
did not include the complex canonical transformations $\exp(i\eps
H_G)$, coming from the averaging procedures along the flows of
$H_p$ and $H_{\langle{s}\rangle}$.

We summarize the discussion above in the following proposition.

\begin{prop}
Assume that the subprincipal symbol of $P_{\eps=0}$ vanishes, and
that the global assumption {\rm (\ref{H})} holds. Let $\eps$
satisfy $h^{1/2}\ll \eps={\cal O}(h^{\delta})$, $\delta>0$. Then
there exists an IR-manifold $\Lambda\subset T^*\widetilde{M}$ and
a smooth Lagrangian torus $\Lambda_{0,0}\subset \Lambda$ such that
when $\rho\in \Lambda$ is away from a small enough \neigh{} of
$\Lambda_{0,0}$ in $\Lambda$ and $\abs{\Re P_{\eps}(\rho,h)}\leq
1/C$ for $C>0$ large enough, it is true that \begeq \label{0.46}
\abs{\frac{\Im
P_{\eps}(\rho,h)}{\eps^3}-\langle\langle{t}\rangle\rangle_{\infty}}\geq
\frac{1}{{\cal O}(1)}.
\endeq
The \neigh{} can be taken as small as we wish provided that $T$ in
{\rm (\ref{H})} and the implicit constant in {\rm (\ref{0.46})}
are chosen large enough, and $h$ is sufficiently small. The
manifold $\Lambda$ is an $(\eps+h/\eps)$-perturbation of $T^*M$
and agrees with this set away from $p^{-1}(0)\cap T^*M$.
Microlocally, near $\Lambda_{0,0}$, the operator
$P_{\eps}:H(\Lambda)\rightarrow H(\Lambda)$ is conjugated by means
of an elliptic uniformly bounded Fourier integral operator to the
translation invariant operator $\widehat{P}$ on
$L^2_{\theta}(\T^2)$.
\end{prop}

Using this proposition and the methods of section 6
of~\cite{HiSj}, we get the final spectral result of this case.

\begin{theo}
We keep the general assumptions of Theorem {\rm 4.2} and consider
the range $h^{1/2}\ll \eps={\cal O}(h^{\delta})$, $\delta>0$.
Assume that the averaged corrections $\langle{s}\rangle$ and
$\langle{t}\rangle$ are real on the real domain, and that the
differentials of $p$ and $\langle{s}\rangle$ are linearly
independent along the set
$$
\Lambda_{0,0}: p=0,\,\langle{s}\rangle=0,
$$
which is assumed to be connected. We assume that $T(0)$ is the
minimal period of every closed $H_p$--trajectory in
$\Lambda_{0,0}$, and we write $S\in \real^2$ and $\alpha\in \z^2$
to denote the actions and Maslov indices along the standard
fundamental cycles in $\Lambda_{0,0}$. Introduce
$\langle\langle{t}\rangle\rangle_{\infty}$ as the mean value of
$\langle{t}\rangle$ over $\Lambda_{0,0}$ and make the global
assumption {\rm (\ref{H})}. Let $C>0$ be large enough. Then the
eigenvalues of $P_{\eps}$ in the rectangle
$$
\abs{\Re z}\leq \frac{1}{C},\quad \abs{\Im z-\eps^3 \langle\langle{t}\rangle\rangle_{\infty}}\leq \frac{\eps^3}{C}
$$
are given, modulo ${\cal O}(h^{\infty})$, by the quasi-eigenvalues
$$
\widehat{P}\left(h(k-\frac{\alpha}{4})-\frac{S}{2\pi},\eps,\frac{h}{\eps};h\right),\quad
k\in \z^2.
$$
Here $\widehat{P}(\xi,\eps,\frac{h}{\eps};h)$ is holomorphic in
$\xi\in {\rm neigh}(0,\comp^2)$, smooth in $\eps,\frac{h}{\eps}\in
{\rm neigh}(0,\real)$, and has an expansion as $h\rightarrow 0$,
$$
\widehat{P}(\xi,\eps,\frac{h}{\eps};h)=p(\xi_1)+\eps^2\left(r_0(\xi,\eps,\frac{h}{\eps})+hr_1(\xi,\eps,\frac{h}{\eps})+\ldots\right).
$$
We have
$$
r_0(\xi,\eps,\frac{h}{\eps})=\langle{s}\rangle(\xi)+{\cal
O}\left(\eps+\frac{h}{\eps}\right),\quad r_j={\cal O}(1),\quad
j\geq 1.
$$
\end{theo}

\medskip \noindent 3. {\it The case when $\eps\sim h^{1/2}$}.

Assume, as we may, that we are now in a situation when the first
two averaging procedures have been carried out, so that we are
dealing with the operator $P_{\eps}$ of principal symbol \begeq
\label{0.47} p_{0,\eps}=p+\eps^2 \langle{s}\rangle+i\eps^3
\langle\langle{t}\rangle\rangle_T+{\cal O}_T(\eps^4),
\endeq
and subprincipal symbol
$$
p_{1,\eps}=\eps q_1+{\cal O}_T(\eps^2),
$$
where $q_1$ does not depend on $\eps$ and $T$. Here we may also assume that the subprincipal symbol has already been
averaged along the flow of $p$, to the leading order, so that
$q_1=\langle{q_1}\rangle$. Since in the relevant parameter range, we
have $\eps^3\sim \eps h$, we shall have to consider not only long-time
averages $\langle\langle{t}\rangle\rangle_T$, but we shall also
average $\langle{q_1}\rangle$ along the
$H_{\langle{s}\rangle}$--flow. Let therefore $G$ solve
$$
H_{\langle{s}\rangle}G=\langle{q_1}\rangle-\langle\langle{q_1}\rangle\rangle_T,
$$
where
$$
\langle\langle{q_1}\rangle\rangle_T=\frac{1}{T}\int_0^T
\langle{q_1}\rangle\circ \exp(uH_{\langle{s}\rangle})\,du,
$$
and compose $P_{\eps}$ with $\exp(\frac{h}{\eps}H_G)$. Then, using
also the fact that $G$ Poisson commutes with $p$, and implementing
$\exp(\frac{h}{\eps}H_G)$ by means of a Fourier integral operator
with the improved Egorov property, we get a new operator with
principal symbol (\ref{0.47}) and with subprincipal symbol
$$
\eps\langle\langle{q_1}\rangle\rangle_T+{\cal O}_T(\eps^2).
$$

In what follows, we assume that $(0,0)$ is a regular value of the
mapping $T^*M\ni \rho\mapsto(p(\rho),\langle{s}\rangle(\rho))\in
\real^2$, with the connected pre-image, and we shall work
microlocally near the Lagrangian torus
$$
\Lambda_{0,0}: p=0,\,\langle{s}\rangle=0.
$$
As before, we shall also assume that $T(0)$ is the minimal period
of the $H_p$--flow in $\Lambda_{0,0}$.

Introducing the quantities
$\langle\langle{t}\rangle\rangle_{\infty}$ and $
\langle\langle{q_1}\rangle\rangle_{\infty}$, defined as the mean
values over $\Lambda_{0,0}$ of $\langle{t}\rangle$ and
$\langle{q_1}\rangle$, respectively, we recall next that as
$T\rightarrow \infty$,
$$
\langle\langle{t}\rangle\rangle_T=\langle\langle{t}\rangle\rangle_{\infty}+{\cal
O}\left(\frac{1}{T}\right),\quad
\langle\langle{q_1}\rangle\rangle_T=\langle\langle{q_1}\rangle\rangle_{\infty}+{\cal
O}\left(\frac{1}{T}\right),
$$
uniformly along $\Lambda_{0,0}$. Using the notation of case 2, we
introduce now the following global assumption, which is a direct
analogue of the hypothesis (\ref{H}):
\begin{eqnarray}
\label{H1} & & \hbox{For any}\,\,b\in (0,a)\,\,\hbox{there exists}
\\ \nonumber & & T(b)>0\,\,\hbox{and}\,\,C(b)>0\,\,\hbox{such
that}\\ \nonumber & &
\inf_{\Lambda_{[b,a]}}\biggl(\langle{\langle{\Im
q_1}\rangle}\rangle_T(\rho) -\langle\langle{\Im q_1}\rangle\rangle_{\infty}\biggr) \geq \frac{1}{C(b)},\quad T\geq T(b),\\
\nonumber & & \sup_{\Lambda_{[-a,-b]}}\biggl(\langle{\langle{\Im
q_1}\rangle}\rangle_T(\rho)-\langle\langle{\Im
q_1}\rangle\rangle_{\infty}\biggr)\leq -\frac{1}{C(b)},\quad T\geq
T(b).
\end{eqnarray}
In what follows, we shall be working under the assumption obtained
by taking the conjunction of the assumptions (\ref{H}) and
(\ref{H1}).

\vskip 2mm \noindent After a passage to the standard torus $\T^2$
and an elimination of the $x_1$--variable, we obtain an operator
$$
P_{\eps}=\sum_{j=0}^{\infty} h^j p_j(x_2,\xi,\eps,T),
$$
with
$$
p_0=p(\xi_1)+\eps^2 \langle{s}\rangle(\xi)+i\eps^3
\langle\langle{t}\rangle\rangle_T(x_2,\xi)+{\cal O}_T(\eps^4),
$$
and
$$
p_1=\eps\langle\langle{q_1}\rangle\rangle_T(x_2,\xi)+{\cal
O}_T(\eps^2).
$$
At this point we may repeat the argument of case 2 to conclude that
there exists $B_0=b_0(x_2,hD_x,\eps,\frac{h}{\eps},T)$ such that
$$
b_0={\cal O}\left(\frac{\eps}{T}\right)+{\cal
O}\left(\frac{h}{\eps T}\right)+ {\cal
O}_T\left((\eps,\frac{h}{\eps})^2\right),
$$
and such that with
$$
r_0(x_2,\eps,\frac{h}{\eps},T)=\langle{s}\rangle(\xi)+i\eps\langle\langle{t}\rangle\rangle_T(x_2,\xi)+{\cal
O}_T(\eps^2)+\frac{h}{\eps}\langle\langle{q_1}\rangle\rangle_T(x_2,\xi)+\frac{h^2}{\eps^2}p_2(x_2,\xi,\eps,T),
$$
the symbol
$$
r_0\circ \exp(H_{b_0})=\sum_{k=0}^{\infty}\frac{H_{b_0}^k}{k!}r_0
$$
is independent of $x_2$. The construction of the operator $B_1$ is
the same as before, and we then introduce a globally defined
IR-manifold $\widetilde{\Lambda}\subset T^*\widetilde{\T}^2$ which
agrees with $\exp(H_{b_0})\left({\rm neigh}(\xi=0,
T^*\T^2)\right)$ near $\xi=0$, is equal to $T^*\T^2$ further out,
and with the property that along $\widetilde{\Lambda}$ we have
$\Im \xi_1=0$ and
$$
\Im x,\,\,\Im \xi_2={\cal O}\left(\frac{\eps}{T}\right)+{\cal
O}\left(\frac{h}{\eps T}\right)+ {\cal
O}_T\left((\eps,\frac{h}{\eps})^2\right).
$$
The action of $P_{\eps}$ on the corresponding Hilbert space
$H(\widetilde{\Lambda})$ is microlocally near $\xi=0$ unitarily
equivalent to the conjugated operator acting on the $L^2$-space of
Floquet periodic functions. The Weyl symbol of the conjugated
operator is independent of $x$ and has the form
$$
p(\xi_1)+\eps^2\left(r_0(\xi,\eps,\frac{h}{\eps},T)+h
r_1(\xi,\eps,\frac{h}{\eps},T)+\ldots\right),
$$
with
$$
r_0=\langle{s}\rangle(\xi)+i\eps\langle\langle\langle{t}\rangle\rangle_T\rangle(\xi)+
\frac{h}{\eps}\langle\langle\langle{q_1}\rangle\rangle_T\rangle(\xi)+{\cal
O}_T\left(\eps^2+\frac{h^2}{\eps^2}\right),
$$
and
$$
r_j={\cal O}_T(1),\quad j\geq 1.
$$
The imaginary part of $P_{\eps}$ along $\widetilde{\Lambda}$ is given by
$$
\eps^3\langle\langle{t}\rangle\rangle_T(\Re x_2,\Re \xi)+{\cal
O}\left(\frac{\eps^3}{T}\right)+h\eps\Im
\langle\langle{q_1}\rangle\rangle_T(\Re x_2,\Re \xi)+{\cal
O}\left(\frac{h\eps}{T}\right)+{\cal O}_T(\eps^4).
$$
Here $\eps\sim h^{1/2}$. It follows from the assumptions (\ref{H})
and (\ref{H1}) that choosing $T$ large enough but fixed, we can
achieve that away from a small \neigh{} of the torus, the modulus
of $\Im
P_{\eps}-\eps^3\langle\langle{t}\rangle\rangle_{\infty}-h\eps
\langle\langle{\Im q_1}\rangle\rangle_{\infty}$ is bounded from
below by $(\eps^3+h\eps)/{\cal O}(1)\sim \eps^3/{\cal O}(1)$, when
the attention is restricted to the region where $\Re P_{\eps}$ is
small. At this stage, the situation is completely analogous to the
previously analyzed case 2, and the rest of the argument goes
through as in that case, without any change.

\vskip 4mm \noindent The following is the main spectral result in
case 3.

\begin{theo}
Let us keep the general assumptions of Theorem {\rm 4.2}. In
particular, assume that the subprincipal symbol of $P_{\eps=0}$
vanishes, so that if $p_{1,\eps}$ is the subprincipal symbol of
$P_{\eps}$ then $p_{1,\eps}=\eps q_1+{\cal O}(\eps^2)$, where
$q_1$ does not depend on $\eps$. Let $\eps\sim h^{1/2}$. Assume
next that the corrections $\langle{s}\rangle$ and
$\langle{t}\rangle$ are real on the real domain, and that the
differentials of $p$ and $\langle{s}\rangle$ are linearly
independent along the set
$$
\Lambda_{0,0}: p=0,\langle{s}\rangle=0.
$$
We assume that $\Lambda_{0,0}$ is connected and that $T(0)$ is the
minimal period for the $H_p$-flow in $\Lambda_{0,0}$. As usual, we
shall write $S\in \real^2$ and $\alpha\in \z^2$ to denote the
classical actions and Maslov indices along the fundamental cycles
in $\Lambda_{0,0}$, with the first cycle corresponding to a closed
$H_p$--trajectory. When $\langle\langle{t}\rangle\rangle_{\infty}$
and $\langle\langle{\Im q_1}\rangle\rangle_{\infty}$ are defined
as the mean values of $\langle{t}\rangle$ and $\langle{\Im
q_1}\rangle$ along $\Lambda_{0,0}$, respectively, we make the
global assumptions {\rm (\ref{H})} and {\rm (\ref{H1})}. Let $C>0$
be sufficiently large. Then the eigenvalues of $P_{\eps}$ in the
domain
$$
\abs{\Re z}<\frac{1}{C},\quad \abs{\Im
z-\eps^3\langle\langle{t}\rangle\rangle_{\infty}-h\eps\langle\langle{\Im
q_1}\rangle\rangle_{\infty}}<\frac{\eps^3}{C}
$$
are given, modulo ${\cal O}(h^{\infty})$, by the formal
quasi-eigenvalues associated with the Lagrangian torus
$\Lambda_{0,0}$,
$$
p\left(h\left(k_1-\frac{\alpha_1}{4}\right)-\frac{S_1}{2\pi}\right)+\eps^2\sum_{j=0}^{\infty}
h^j r_j
\left(h\left(k-\frac{\alpha}{4}\right)-\frac{S}{2\pi},\eps,\frac{h}{\eps}\right),\quad
k=(k_1,k_2)\in \z^2.
$$
\end{theo}

\section{Barrier top resonances: the case of $1:1$ resonance}\label{section5}
\setcounter{equation}{0}
Consider
\begeq
\label{5.1}
P=-h^2\Delta+V(x),\quad p(x,\xi)=\xi^2+V(x),\quad (x,\xi)\in
T^*\real^2,
\endeq
where $V$ is an analytic potential, satisfying the same general
assumptions as in section 7 of~\cite{HiSj}, which allow us to
define the resonances of $P$ in a fixed sector in the fourth
quadrant. As in~\cite{HiSj}, we assume that $V(0)=E_0>0$,
$V'(0)=0$, $V''(0)<0$, and that $(0,0)$ is the only trapped
$H_p$--trajectory in $p^{-1}(E_0)\cap \real^4$. The Taylor
expansion of $p(x,\xi)$ in suitable linear symplectic coordinates
has the form \begeq \label{5.2} p(x,\xi)-E_0=\sum_{j=1}^2
\frac{\lambda_j}{2}(\xi_j^2-x_j^2)+p_3(x)+p_4(x)+\ldots,\quad
(x,\xi)\rightarrow 0.
\endeq
Here $\lambda_j>0$ and $p_j(x)$ is a homogeneous polynomial of
degree $j\geq 3$. We are interested in resonances of $P$ near
$E_0$, and from~\cite{HiSj} we recall that the study of such
resonances can be reduced to an eigenvalue problem for $P-E_0$ near $0$, after the complex scaling given by
$x=e^{i\pi/4}\widetilde{x}$, $\xi=e^{-i\pi/4}\widetilde{\xi}$,
$(\widetilde{x},\widetilde{\xi})\in T^*\real^2$. Performing the
scaling and dropping the tildes from the notation, we get a new
operator with the leading symbol
$$
\frac{1}{i}\left(p_2(x,\xi)+i e^{3\pi i/4}p_3(x)-ip_4(x)+\ldots
\right)=:\frac{1}{i}q(x,\xi),\quad (x,\xi)\rightarrow 0,
$$
with \begeq \label{5.3} p_2(x,\xi)=\sum_{j=1}^2
\frac{\lambda_j}{2}\left(x_j^2+\xi_j^2\right),
\endeq
and with the vanishing subprincipal symbol.

\vskip 2mm
We shall be interested in eigenvalues $E$ of the
operator
$$
Q(x,hD_x;h)=q(x,hD_x)+{\cal O}(h^2),
$$
with $\abs{E}\sim \eps^2$, $h^{\delta}<\eps\ll 1$, $0<\delta<1/2$.
After a rescaling $x=\eps y$, we get
$$
\frac{1}{\eps^2}Q(x,hD_x;h)=\frac{1}{\eps^2}Q(\eps(y,\widetilde{h}D_y);h),\quad
\widetilde{h}=\frac{h}{\eps^2}\ll 1,
$$
with the corresponding symbol
$$
\frac{1}{\eps^2}Q(\eps(y,\eta))\sim
\frac{1}{\eps^2}q(\eps(y,\eta))+\eps^2 \widetilde{h}^2
q^{(2)}(\eps(y,\eta))+\ldots.
$$
The leading symbol becomes
$$
\frac{1}{\eps^2}q(\eps(y,\eta))=p_2(y,\eta)+i\eps e^{3\pi
i/4}p_3(y)-i\eps^2 p_4(y)+\ldots,
$$
to be considered in a region where $\abs{(y,\eta)}\sim 1$, where
the corresponding eigenfunctions are concentrated.

A straightforward application of Theorem 3.1 together with the
scaling reductions above gives the following result.
\begin{theo}
Assume that the principal symbol $p(x,\xi)$ of the operator
{\rm (\ref{5.1})} has an expansion {\rm (\ref{5.2})} with
\begeq
\label{5.4}
\lambda\cdot k=0,\quad \wrtext{for some}\quad 0\neq k\in \z^2,
\endeq
and assume that the average of $p_3$ along the closed
$H_{p_2}$--trajectories, $\langle{p_3}\rangle$, vanishes
identically. Introduce the function
$$
s=-p_4+\frac{1}{2T}\int_0^T t\{p_3\circ \exp(tH_{p_2}),p_3\}\,dt,
$$
where $T>0$ is the period of the $H_{p_2}$--flow on $p_2^{-1}(1)$,
and assume that $\langle{s}\rangle$ is not identically zero. Then
the resonances $E$ of the operator {\rm (\ref{5.1})} in the domain
\begeq \label{5.5} \{z\in \comp; h^{2/3}\ll \abs{z-E_0}={\cal
O}(1)h^{\delta}\}\backslash \bigcup \{z\in \comp; \abs{\Re
z-E_0-A\abs{\Im z}^2}<\eta \abs{\Im z}^2\},
\endeq
where $\delta$, $\eta>0$ are arbitrary but fixed, are given by
\begeq \label{5.5.1} E=E_0-i\eps^2
\widehat{P}\left(\widetilde{h}\left(k-\frac{\alpha}{4}\right)-\frac{S}{2\pi},\eps,\frac{\widetilde{h}}{\eps};\widetilde{h}\right)+{\cal
O}(h^{\infty}),\quad \widetilde{h}=\frac{h}{\eps^2},\,\,k\in \z^2.
\endeq
Here we choose $\eps>0$ with $\abs{E-E_0}\sim \eps^2$ and the
union in {\rm (\ref{5.5})} is taken over the set of critical
values of $\langle{s}\rangle$ restricted to $p_2^{-1}(1)$, with
$A$ varying over this set.

The function
$\widehat{P}(\xi,\eps,\frac{\widetilde{h}}{\eps};\widetilde{h})$
in {\rm(\ref{5.5.1})} has an asymptotic expansion as
$\widetilde{h}\rightarrow 0$,
$$
\widehat{P}(\xi,\eps,\frac{\widetilde{h}}{\eps};\widetilde{h})\sim
p(\xi_1)+\eps^2\sum_{j=0}^{\infty}\widetilde{h}^j
r_j\left(\xi,\eps,\frac{\widetilde{h}}{\eps}\right),
$$
where
$$
r_0=i\langle{s}\rangle+{\cal
O}\left(\eps+\frac{\widetilde{h}}{\eps}\right),\,\,r_j={\cal
O}\left(\eps+\frac{\widetilde{h}}{\eps}\right),\,\,j\geq 1.
$$
\end{theo}

In section 7 of~\cite{HiSj} we considered the case of the resonant
frequencies $(\lambda_1,\lambda_2)=(1,2)$, and we have also
remarked there that in the case of the $1:1$ resonance,
$(\lambda_1,\lambda_2)={\rm Const}(1,1)$, it is true that
$\langle{p_3}\rangle\equiv 0$, for any cubic polynomial $p_3$. We
shall now illustrate Theorem 5.1 by discussing this explicit
example, where we shall take $(\lambda_1,\lambda_2)=(1,1)$. In
doing so, we shall also assume for simplicity that $p_4\equiv 0$.

\noindent
\vskip 2mm
It will be convenient to work in the symplectic coordinates $(y,\eta)$ given by
$$
y=\frac{1}{\sqrt{2}}(x-i\xi),\quad \eta=\frac{1}{i\sqrt{2}}(x+i\xi).
$$
In these coordinates we have $p_2=\sum_{j=1}^2 i\lambda_j y_j \eta_j$,
and the flow is given by
$$
\exp(tH_{p_2})(y,\eta)=(e^{it\lambda_1}y_1,
e^{it\lambda_2}y_2,e^{-it\lambda_1}\eta_1,e^{-it\lambda_2}\eta_2).
$$
When $\abs{\alpha}=3$, we write
$$
p_3(x)=x^{\alpha}=\sum_{0\leq k\leq \alpha}a_{k\alpha} y^k \eta^{\alpha-k},
$$
where
$$
a_{k\alpha}=\frac{i^{\abs{\alpha-k}}}{2^{\abs{\alpha}/2}}
\pmatrix{\alpha \cr k}.
$$
When computing
$$
G_0=\frac{1}{2\pi}\int_0^{2\pi} t p_3\circ \exp(tH_p)\,dt,
$$
we use that
$$
{1\over 2\pi}\int_0^{2\pi} t e^{it\lambda \cdot (2k-\alpha)}dt
=\cases{\pi \hbox{ if }\lambda \cdot (2k-\alpha)=0,\cr \frac{1}{i(\lambda\cdot(2k-\alpha))}\hbox{ otherwise,}}
$$
and obtain
$$
G_0=\sum_{0\leq k\leq \alpha \atop {2k=\alpha+nk_0}}\pi a_{k\alpha}y^k
\eta^{\alpha-k}+\sum_{0\leq k\leq \alpha \atop{2k-\alpha\neq nk_0}}\frac{a_{k\alpha}}{i\lambda\cdot(2k-\alpha)}
y^k\eta^{\alpha-k},\quad n\in \z,
$$
which we write as
$$
G_0=\sum_{0\leq k\leq \alpha}g_{k\alpha}y^k \eta^{\alpha-k}.
$$
Here $k_0\in \z^2$ satisfies (\ref{5.4}) and has the minimal norm and positive first component.
We are then interested in computing $s=(1/2)H_{G_0}p_3$, which is
equal to
\begin{eqnarray}
\label{par}
& & \half H_{G_0}p_3=\half \sum_{0\leq k\leq \alpha \atop{0\leq \widetilde{k}\leq
\alpha}}g_{k\alpha} a_{\widetilde{k}\alpha}
\{y^{k}\eta^{\alpha-k},y^{\widetilde{k}},\eta^{\alpha-\widetilde{k}}\}
\\ \nonumber
& = & \half \sum_{0\leq k\leq \alpha \atop{0\leq \widetilde{k}\leq \alpha}}g_{k\alpha}
a_{\widetilde{k}\alpha} \sum_{j=1}^2
\sigma(k_j,\alpha_j-k_j;\widetilde{k}_j,\alpha_j-\widetilde{k}_j)y^{k+\widetilde{k}-e_j}\eta^{2\alpha-k-\widetilde{k}-e_j}.
\end{eqnarray}
Here $\sigma$ is the symplectic form on $\real^4$ and $e_1=(1,0)$,
$e_2=(0,1)$. Now when computing the flow average of (\ref{par}),
we notice that the only non-vanishing contribution will come from the
terms for which $2(k+\widetilde{k})=2\alpha+nk_0$ for some $n\in
\z$. Therefore,
\begeq
\label{5.6.5}
\langle{s}\rangle=\half \langle{H_{G_0}p_3}\rangle=\half \sum_{0\leq k\leq \alpha \atop{0\leq
\widetilde{k}\leq \alpha \atop{2(k+\widetilde{k})=2\alpha+nk_0}}}
g_{k\alpha} a_{\widetilde{k}\alpha} \sum_{j=1}^2
\alpha_j(k_j-\widetilde{k}_j)y^{k+\widetilde{k}-e_j}\eta^{2\alpha-k-\widetilde{k}-e_j}.
\endeq
Here we have also used that
$$
\sigma(k_j,\alpha_j-k_j;\widetilde{k}_j,\alpha_j-\widetilde{k}_j)=
\sigma(k_j,\alpha_j;\widetilde{k}_j,\alpha_j)=\alpha_j(k_j-\widetilde{k}_j).
$$

We shall now consider the case when $(\lambda_1,\lambda_2)=(1,1)$ and
$\alpha=(\abs{\alpha},0)$ so that
$$
p_3(x)=x_1^3=2^{-3/2}\left(y^3_1+3iy^2_1\eta_1-3y_1\eta_1^2-i\eta_1^3\right).
$$
Then
$$
a_{(3,0);(3,0)}=2^{-3/2},\quad a_{(2,0);(3,0)}=2^{-3/2}3i,
$$
$$
a_{(1,0);(3,0)}=-2^{-3/2}3,\quad a_{(0,0);(3,0)}=-2^{-3/2}i.
$$
and since $\lambda\cdot(2k-\alpha)=2\abs{k}-\abs{\alpha}\neq 0$, we
see that
$$
g_{k\alpha}=\frac{a_{k\alpha}}{i(2k_1-3)}.
$$
Now $k_0=(1,-1)$, and it follows that in (\ref{5.6.5}) we have $n=0$
and $k+\widetilde{k}=\alpha$. We get
$$
\langle{H_{G_0}p_3}\rangle=\frac{3 y_1^2 \eta_1^2}{i}\sum_{0\leq k\leq
\alpha}a_{k\alpha}a_{\alpha-k,\alpha},
$$
and a straightforward computation shows that
$$
\langle{s}\rangle=-\frac{15}{4}y_1^2\eta_1^2=\frac{15}{4}\left(\frac{x_1^2+\xi_1^2}{2}\right)^2.
$$
We then see that the differential of this function restricted to the sphere
$$
p_2(x,\xi)=\frac{x_1^2+\xi_1^2}{2}+\frac{x_2^2+\xi_2^2}{2}=1
$$
vanishes along two closed $H_{p_2}$--trajectories given by
\begeq
\label{5.6}
x_1=\xi_1=0,\quad \frac{x_2^2+\xi_2^2}{2}=1,
\endeq
and \begeq \label{5.7} x_2=\xi_2=0,\quad
\frac{x_1^2+\xi_1^2}{2}=1,
\endeq
with the corresponding critical values $0$ and $15/4$. Notice also
that the critical trajectory (\ref{5.6}) is degenerate in the
sense that the transversal Hessian of $\langle{s}\rangle$ is
degenerate. An application of Theorem 5.1 then gives a description
of all the resonances of $P$ in an energy shell of the form
$h^{2/3}\ll \abs{z-E_0}={\cal O}(h^{\delta})$, $\delta>0$, after
we have deleted arbitrarily small parabolic neighborhoods of the
curves
$$\{z\in\comp; \Re z=E_0+A\abs{\Im z}^2\},
$$
where $A\in \{0,15/4\}$.

\section{Complex perturbations on the 2-sphere}\label{section6}
\setcounter{equation}{0}

On the sphere $S^2=\{x\in \real^3; x_1^2+x_2^2+x_3^2=1\} \subset
\real^3$, equipped with the standard metric, let us consider an
operator of the form $P_{\eps}=-h^2\Delta+i\eps q(x)$, where
$q(x)$ is a real-valued analytic function. We introduce the
leading symbol of $P_{\eps=0}$, $p$, on $T^*S^2$, and recall that
the geodesic flow, which can be identified with the Hamilton flow
of $p$, is periodic on each energy surface $p^{-1}(E)$, $E>0$,
with the minimal period \begeq \label{6.1}
T(E)=\frac{\pi}{\sqrt{E}}.
\endeq
We shall now describe the $H_p$--flow on $T^*S^2$ in a more
explicit way. In doing so, let us remark that the Euclidean metric
on $\real^3$ allows us to embed $T^*S^2$ symplectically into
$T^*\real^3$, $T^*S^2\simeq TS^2\hookrightarrow T\real^3 \simeq
T^*\real^3$, so that
$$
T^*S^2\simeq \Sigma :=\{(x,\xi)\in T^*\real^3; h_1(x,\xi)=0,
h_2(x,\xi)=0\},
$$
where $h_1=x^2-1$ and $h_2=x\cdot\xi$, and  where we use the
Euclidean scalar product and the norm. In what follows we shall
write $\{\cdot,\cdot\}$ to denote the Poisson bracket on
$T^*\real^3$ and the Poisson bracket on $\Sigma$ will be denoted by
$\{\cdot,\cdot\}_{\Sigma}$.

$\Sigma$ is a symplectic submanifold of $T^*\real^3$, since
$\{h_1,h_2\}=-2$ there. The restriction to $\Sigma$ of the
Liouville form $\xi\cdot dx$ can be identified with the
corresponding Liouville form on $T^*S^2$. We then have the
corresponding identification of the symplectic forms. Furthermore,
$p=\xi^2$, when viewed as a function on $\Sigma$, and the
corresponding Hamilton field satisfies $H_p^{\Sigma} \equiv
H_p\,\,{\rm mod}\,\, T\Sigma^{\sigma}$, where $H_p$ is the
Hamilton field on $T^*\real^3$ and the exponent $\sigma$ indicates
that we take the symplectic orthogonal. It follows that
$H_p^{\Sigma}=H_{\widetilde{p}}$, where $\widetilde{p}\in
C^{\infty}(\real^3)$ is such that $p=\widetilde{p}$ along $\Sigma$
and $H_{\widetilde{p}}$ is tangent to $\Sigma$. Writing
$\widetilde{p}=p+ah_1+bh_2$, where $a$ and $b$ are chosen so that
$\{\widetilde{p},h_j\}=0$ along $\Sigma$, $j=1,2$, we find after a
simple computation that $\widetilde{p}=x^2\xi^2$. Restricting to
$x^2=1$, we get
$$
\frac{1}{2}H_{x^2\xi^2}=\xi\cdot \partial_x-\xi^2 x\cdot
\partial_{\xi}.
$$
This implies that the map $\exp(tH_p^{\Sigma}/2): T^*S^2\backslash 0
\rightarrow T^*S^2\backslash 0 $ is given by \begeq \label{6.2}
(x,\xi)\mapsto
(\cos(\abs{\xi}t)x+\sin(\abs{\xi}t)\frac{\xi}{\abs{\xi}},
-\abs{\xi}\sin(\abs{\xi}t)x+\cos(\abs{\xi}t)\xi).
\endeq

\noindent {\it Example.} When $q(x)=x_1x_2$, we shall compute the
averaged perturbation \begeq \label{6.3} \langle{q}\rangle =
\frac{1}{2\pi}\int_0^{2\pi} q\circ \exp(tH_p^{\Sigma}/2)\,dt\quad
\wrtext{on}\quad p^{-1}(1),
\endeq
and using (\ref{6.2}) we obtain
 \begeq \label{6.4}
\langle{q}\rangle(x,\xi)=\frac{x_1 x_2+\xi_1 \xi_2}{2}.
\endeq
Now $\langle{q}\rangle$ is invariant under the
$H_p^{\Sigma}$-flow, and should therefore be viewed as a function
on the reduced space of oriented closed orbits in $p^{-1}(1)$,
\begeq \label{6.5} {\cal O}:= p^{-1}(1)/\exp(\real H_p^{\Sigma}).
\endeq
Here ${\cal O}$ is a two-dimensional compact symplectic manifold,
which consists of all oriented great circles of $S^2$, and we can
identify ${\cal O}$ with $S^2$ by associating the unit vector
$y=x\times \xi$ to the $H_p^{\Sigma}$--trajectory in $p^{-1}(1)$,
given by (\ref{6.2}). In order to express the flow average
$\langle{q}\rangle$ in terms of $y$, we notice that we can choose
$(x,\xi)$ corresponding to $y$ with \begeq \label{6.6.1}
x=\frac{(-y_2,y_1,0)}{(y_1^2+y_2^2)^{1/2}},
\endeq
and \begeq \label{6.6.2} \xi=y\times
x=\frac{(-y_1y_3,-y_2y_3,y_1^2+y_2^2)}{(y_1^2+y_2^2)^{1/2}}.
\endeq
It follows that \begeq \label{6.7} \langle{q}\rangle(y)=-\frac{y_1
y_2}{2}.
\endeq

\vskip 2mm \noindent A straightforward computation shows that
$\langle{q}\rangle$ is a Morse function on $S^2$, with 6 critical
points. We get two non-degenerate maxima at $\pm (1/\sqrt{2},
-1/\sqrt{2},0)$, with the corresponding critical value $1/4$, two
non-degenerate minima at $\pm (1/\sqrt{2},1/\sqrt{2},0)$ with the
critical value $-1/4$, and two non-degenerate saddle points at the
poles $\pm (0,0,1)$, with the corresponding critical value 0.

\Remark. The considered example comes essentially from~\cite{Gr}.
Notice also that the expression (\ref{6.7}) is in agreement with
the general fact, established in~\cite{Gu1}, that the Radon
transform $q\mapsto \langle{q}\rangle$ on $S^2$ maps the space of
the restrictions of the homogeneous harmonic polynomials of a
fixed degree into itself, and is a multiple of the identity when
restricted to this space. A simple proof of this fact
following~\cite{Gu1}, is obtained if we observe that the Radon
transform commutes with the SO(3) action on $S^2$, and apply
Schur's lemma. Let us also recall from~\cite{He} and~\cite{Gu1}
that the kernel of the Radon transformation, viewed as a mapping
on $C^{\infty}(S^2)$, consists precisely of all odd functions on
$S^2$, and that it becomes bijective when restricted to the space
of smooth even functions. It follows that every smooth Morse
function in the range of the Radon transform has at least two
saddle points.

\vskip 4mm \noindent As a preparation for further considerations,
we shall now derive an expression for the Poisson bracket of two
functions $f, g\in C^{\infty}(T^*\real^3)$, viewed as functions on
the symplectic manifold $\Sigma$. In doing so, we have to compute the
Hamilton field of the restriction of $f$ to $\Sigma$, and repeating the
previous arguments, we find that $H_f^{\Sigma}=H_{\widetilde{f}}$, where
$$
\widetilde{f}=f+\frac{\{f,h_2\}}{\{h_2,h_1\}}h_1+\frac{\{f,h_1\}}{\{h_1,h_2\}}h_2.
$$
Therefore,
\begeq \label{6.8}
\{f,g\}_{\Sigma}=\{f,g\}+\frac{1}{2}\left(\{f,h_2\}\{h_1,g\}-\{f,h_1\}\{h_2,g\}\right).
\endeq

We summarize the discussion above in the following essentially
well-known pro\-po\-si\-tion---see also~\cite{DLT} for a more
general discussion of Hamilton mechanics with constraints.
\begin{prop}
Let
$$
\Sigma=\{(x,\xi)\in T^*\real^3; h_1(x,\xi)=0, h_2(x,\xi)=0\},\quad
h_1=x^2-1,\quad h_2=x\cdot \xi.
$$
Then $\Sigma$ is a symplectic submanifold of $T^*\real^3$ and the
Hamilton vector field of $p=\xi^2$, viewed as a function on $\Sigma$,
is given by $H_{x^2\xi^2}$. The Poisson bracket of the
restrictions of $f,g\in C^{\infty}(T^*\real^3)$ to $\Sigma$ is given by
{\rm (\ref{6.8})}.
\end{prop}

We now consider the case when the perturbation $q$ is an odd
function in $\real^3$. Then $\langle{q}\rangle\equiv 0$, and
assuming that $q$ is an odd monomial, we shall derive an
expression for the second averaged correction $\langle{s}\rangle$,
given by (\ref{2.6}). In doing so, when $(x,\xi)\in
T^*S^2\backslash 0 $, we introduce $z=x+i\xi/\abs{\xi}\in
\comp^3$. It follows then from (\ref{6.2}) that along the
$H_{p/2}^{\Sigma}$-trajectory we have \begeq \label{6.9}
 z(t)=e^{-i\abs{\xi}t}z(0).
\endeq
Therefore, with $q(x)=x^{\alpha}$, $\abs{\alpha}$ odd, we get
$$
q(\exp(t/2H_p^{\Sigma})(x,\xi))=\left(\frac{ze^{-i\abs{\xi}t}+\overline{z}e^{i\abs{\xi}t}}{2}\right)^{\alpha}=\frac{1}{2^{\abs{\alpha}}}\sum_{0\leq
\beta\leq \alpha}\pmatrix{\alpha \cr
\beta}z^{\beta}\overline{z}^{\alpha-\beta}e^{it\abs{\xi}(\abs{\alpha}-2\abs{\beta})}.
$$
We now have to compute the function $G_0$, given by (\ref{2.2}),
$$
G_0=\frac{1}{T(E)}\int_0^{T(E)}t q\circ\exp(tH_p)\,dt \quad
\wrtext{on}\quad p^{-1}(E),
$$
and we immediately see that \begeq \label{6.10}
G_0=\frac{1}{2^{\abs{\alpha}+1}\abs{\xi}}\sum_{0\leq \beta\leq
\alpha}\pmatrix{\alpha \cr
\beta}\frac{z^{\beta}\overline{z}^{\alpha-\beta}}{i(\abs{\alpha}-2\abs{\beta})}.
\endeq
We shall next consider the Poisson bracket of $G_0$ and
$$
q(x)=x^{\alpha}=\frac{1}{2^{\abs{\alpha}}}\sum_{0\leq \gamma\leq
\alpha}\pmatrix{\alpha \cr
\gamma}z^{\gamma}\overline{z}^{\alpha-\gamma},
$$
so that
$$
\{G_0,q\}_{\Sigma}=\frac{1}{2^{2\abs{\alpha}+1}}\sum_{0\leq \beta\leq
\alpha\atop{0\leq \gamma\leq \alpha}} \pmatrix{\alpha \cr
\beta}\pmatrix{\alpha \cr
\gamma}\frac{1}{i(\abs{\alpha}-2\abs{\beta})}
\biggl\{\frac{z^{\beta}\overline{z}^{\alpha-\beta}}{\abs{\xi}},z^{\gamma}\overline{z}^{\alpha-\gamma}\biggr\}_{\Sigma}.
$$
When computing the flow average of this expression, we see using
(\ref{6.9}) and the fact that $\abs{\xi}$ is constant along the
$H_p^{\Sigma}$--flow, that the only non-zero contributions to the average
come from the terms for which
$\abs{\alpha}=\abs{\beta}+\abs{\gamma}$. We get
$$
\langle{\{G_0,q\}_{\Sigma}}\rangle=\frac{1}{2^{2\abs{\alpha}+1}}\sum_{0\leq
\beta\leq\alpha\atop{0\leq \gamma\leq
\alpha\atop{\abs{\alpha}=\abs{\beta}+\abs{\gamma}}}}\pmatrix{\alpha
\cr \beta}\pmatrix{\alpha \cr
\gamma}\frac{1}{i(\abs{\alpha}-2\abs{\beta})}\biggl\{\frac{z^{\beta}\overline{z}^{\alpha-\beta}}{\abs{\xi}},
z^{\gamma}\overline{z}^{\alpha-\gamma}\biggr\}_{\Sigma}.
$$

We shall now illustrate the preceding discussion by an explicit
example when $\alpha=(1,0,0)$, $q(x)=x_1$. Using (\ref{6.10}) we
see that
$$
G_0(x,\xi)=-\frac{\xi_1}{2\xi^2},
$$
and we may then check directly that $H_p^{\Sigma} G_0=q$. A direct
computation using (\ref{6.8}) shows next that
$$
\{G_0,q\}_{\Sigma}=\frac{1}{2}\left(\frac{x_1^2-1}{\xi^2}+\frac{2\xi_1^2}{\xi^4}\right).
$$
When computing the flow average of $\{G_0,q\}_{\Sigma}$ restricted to
$p^{-1}(1)$, we notice that $x_1^2+\xi_1^2$ is flow invariant,
while using (\ref{6.2}) we find that the flow average of $\xi_1^2$
is $(x_1^2+\xi_1^2)/2$. Therefore,
$$
\langle{s}\rangle=-\frac{1}{2}\langle{\{G_0,q\}_{\Sigma}}\rangle=\frac{1}{4}-\frac{3(x_1^2+\xi_1^2)}{8},\quad
\wrtext{on}\quad p^{-1}(1),
$$
and using (\ref{6.6.1}) and (\ref{6.6.2}), we find that when
viewed as a function on ${\cal O}\simeq S^2$, $\langle{s}\rangle$
becomes \begeq \label{6.11}
\langle{s}\rangle=\frac{3}{8}y_1^2-\frac{1}{8}.
\endeq
The differential of this function vanishes along the equator
$y_2^2+y_3^2=1$, with the corresponding critical value $-1/8$, and
at the poles $\pm(1,0,0)$ with the critical value $1/4$.

Combining Proposition 2.1 and remark at the end of section 2,
together with Theorem 4.2 and the discussion above, we obtain the
following result.

\begin{theo}
When $q(x)=x_1$ on $\real^3$, let us consider the operator
$P_{\eps}=-h^2\Delta+i\eps q(x)$ acting on $L^2(S^2)$. Assume that
$\eps\ll h^{1/2}$. Then the spectrum of $P_{\eps}$ in ${\rm
neigh}(1,\comp)$ is contained in the union of the rectangles of
the form
$$
h^2\left(k+\frac{1}{2}\right)^2+[-{\cal O}(\eps^2+h^2), {\cal
O}(\eps^2+h^2)]+i[-{\cal O}(\eps h), {\cal O}(\eps h)],\quad k\in
\z.
$$
Let $F_0\in (-1/8, 1/4)$ and restrict the range of $\eps$ further
by assuming that $h\ll \eps$. Then the eigenvalues of $P_{\eps}$
in the set
$$
\abs{\Re z-h^2\left(k+\frac{1}{2}\right)^2-\eps^2 F_0}<
\frac{\eps^2}{{\cal O}(1)},\quad k\in \z,
$$
are given by
$$
\sim h^2\left(k+\frac{1}{2}\right)^2+\eps^2\sum_{j=0}^{\infty} h^j
r_j\left(h\left(k+\frac{1}{2}\right)-1,
h\left(l-\frac{\alpha}{4}\right)-\frac{S}{2\pi},\eps,\frac{h}{\eps}\right),\quad
l\in \z.
$$
Here
$$
r_0\left(\xi,\eps,\frac{h}{\eps}\right)=\langle{s}\rangle(\xi)+{\cal
O}\left(\eps+\frac{h}{\eps}\right),\quad
r_j\left(\xi,\eps,\frac{h}{\eps}\right)={\cal
O}\left(\eps+\frac{h}{\eps}\right),\quad j\geq 1,
$$
are holomorphic in $\xi\in {\rm neigh}(0,\comp^2)$ and smooth in
$\eps,h/\eps\in {\rm neigh}(0,\real)$. The coordinates
$\xi_1=\xi_1(E)$ and $\xi_2=\xi_2(E,F)$ are the normalized actions
of
$$
\Lambda_{E,F}: p=E,\, \langle{s}\rangle=F,
$$
for $E\in {\rm neigh}(1,\real)$, $F\in {\rm neigh}(F_0,\real)$,
given by
$$
\xi_j=\frac{1}{2\pi}\left(\int_{\gamma_j(E,F)}\eta\,dy-\int_{\gamma_j(1,F_0)}\eta
\,dy\right),\quad j=1,2,
$$
with $\gamma_j(E,F)$ being fundamental cycles in $\Lambda_{E,F}$,
such that $\gamma_1(E,F)$ is given by a closed $H_p$--trajectory
of minimal period $T(E)$. We have $\xi_1(E)=E^{1/2}-1$, so that
when expressed in terms of $\xi_1$, $p$ becomes $(\xi_1+1)^2$.
Finally, $\alpha\in \z$ and $S\in \real$ are fixed.
\end{theo}

\Remark. We shall finish this section by pointing out that the
methods of~\cite{HiSj} and of the present paper can be extended to
cover the case when the operator $P_{\eps}$, introduced in section
1, depends holomorphically on a parameter $z\in {\rm
neigh}(0,\comp)$, in such a way that on a symbol level, we have a
holomorphic family of functions $P(x,\xi,\eps,z;h)$, with all the
properties described in section 1 holding uniformly in $z$, and
such that the leading symbol of $P(x,\xi,0,z;h)$ is of the form
$p(x,\xi)-z$, for a real-valued $p(x,\xi)$ with a periodic
Hamilton flow. We then check by inspection that the reductions of
sections 3 and 4 work in the same way, and as a microlocal
Birkhoff normal form for $P_{\eps}$, we obtain a family of
translation invariant operators on $\T^2$ depending smoothly on
$z\in {\rm neigh}(0,\comp)$. The setup of the global Grushin
problem is the same as before---see also~\cite{MeSj}, where the
general holomorphic dependence on the spectral parameter is
assumed from the beginning. The complete asymptotic expansions for
the eigenvalues then come from applying the implicit function
theorem to the "effective Hamiltonian" $E_{-+}$, as in section 6
of~\cite{MeSj}.

These observations are motivated by the problem of studying
asymptotics of eigenfrequencies associated to the damped wave
equation with an analytic damping coefficient, on an analytic Zoll
surface $M$. Indeed, let us recall from~\cite{Sj3} and~\cite{Hi}
that after a semiclassical reduction in the eigenfrequency
equation
$$
\left(-\Delta+2ia(x)\tau-\tau^2\right)u=0, \quad \abs{\tau}\gg 1,
$$
obtained by writing
$$
\tau=\frac{\sqrt{z}}{h},\quad 0<h\ll 1,\quad z\in {\rm
neigh}(1,\comp),
$$
one is led to consider the problem
$$
\left(-h^2\Delta+2ih\sqrt{z}a(x)-z\right)u=0.
$$
Therefore, in this case $P_{\eps}=-h^2\Delta+2i\eps \sqrt{z}a$
with $\eps=h$, and the subprincipal symbol of
$P_{\eps=0}=-h^2\Delta$ vanishes. Without any analyticity
assumptions and without any restriction on the dimension, it was
proved in~\cite{Hi} that the spectrum of $P_{\eps}$ in ${\rm
neigh}(1,\comp)$ has a cluster structure, for any $a\in
C^{\infty}(M;\real)$. If now $M$ is an analytic compact symmetric
surface of rank one which is not $S^2$, then from~\cite{Gu3} we
know that the flow average $\langle{a}\rangle \not \equiv 0$
unless $a$ vanishes identically. The results of~\cite{HiSj},
generalized to cover the case of holomorphic dependence on the
spectral parameter in the non-selfadjoint perturbation, as
indicated above, give then complete asymptotic expansions for the
eigenvalues $z$, in the sub-clusters corresponding to regular
values of $\langle{a}\rangle$, viewed as a function on
$p^{-1}(1)$, and to non-degenerate extreme values of
$\langle{a}\rangle$ on $p^{-1}(1)$. If $M=S^2$ and $a$ is not an
odd function, then the results of~\cite{HiSj} are still
applicable, since, as was remarked before, in this case
$\langle{a}\rangle$ does not vanish identically. Consider finally
the case when $M=S^2$ and $a$ is an analytic odd function. It
follows then that the results of sections 4 and 6 of the present
paper, generalized as described above, give complete asymptotic
expansions for the eigenvalues $z$ in the sub-clusters
corresponding to regular values of the averaged second correction
$\langle{s}\rangle$, after the damping coefficient $a$ has been
multiplied by a sufficiently large but fixed coupling constant.

\end{document}